\magnification 1100
\def\NC{1.1}
\def\PO{1.2}

\def\BS{2.1}
\def\CS{2.2}
\def\exp{2.2.1}
\def\rig{2.3}
\def\lsy{2.3.1}
\def\AP{2.4}
\def\cor{2.4.1}
\def\num{2.4.2}

\def\CI{3.1}
\def\cyc{3.2}
\def\ecc{3.2.1}
\def\max{3.2.2}
\def\rel{3.3}
\def\ctp{3.3.1}
\def\inv{3.3.2}
\def\esp{3.4}
\def\comb{3.4.1}
\def\split{3.4.2}

\def\ele{4.1}
\def\HYP{4.2}
\def\HAT{4.2.1}
\def\GA{4.3}
\def\git{4.3.1}
\def\sta{4.3.2}
\def\TP{4.4}
\def\AMS{4.4.5}


\def\ET{5.1}
\def\USC{5.1.1}
\def\NB{5.2}
\def\van{5.2.1}
\def\dif{5.2.2}

\def\main{6.1.1}

\font\tenmsb=msbm10
\font\sevenmsb=msbm10 at 7pt
\font\fivemsb=msbm10 at 5pt
\newfam\msbfam
\textfont\msbfam=\tenmsb
\scriptfont\msbfam=\sevenmsb
\scriptscriptfont\msbfam=\fivemsb
\def\Bbb#1{{\fam\msbfam\relax#1}}

\def\sqr#1#2{{\vcenter{\vbox{\hrule height.#2pt \hbox{\vrule width.#2pt
height #1pt \kern #1pt \vrule width.#2pt}\hrule height.#2pt}}}}
\def\qed{$\  \sqr74$}

\def\mapdown#1{\Big\downarrow \rlap{$\vcenter{\hbox{$\scriptstyle#1$}}$}}

\def\lmapdown#1{\llap{$\vcenter{\hbox{$\scriptstyle#1$}}\Big\downarrow $}}
\def\mapright#1{\smash{\mathop{\longrightarrow}\limits^{#1}}}

 \def\la{\longrightarrow}

 \def\ni{\noindent}
 \def\cl{\centerline}
\def\df{\noindent {\bf Definition.\  }}
 \def\rk{\noindent {\it Remark.\  }}
 
 \def\pf{\noindent {\it Proof.\  }}

\def\d{\delta }
 
 \def\l{\lambda}
 \def\s{\sigma}
\def\a{\alpha}
\def\b{\beta}
\def\g{\gamma}
\def\t{\theta}
\def\T{\Theta}

\font\gothicf=eufm10
\font\sgothic=eufm7
\font\ssgothic=eufm5
\textfont5=\gothicf
\scriptfont5=\sgothic
\scriptscriptfont5=\ssgothic

\def\C{{\Bbb C}}

\def\P{{\Bbb P}}
\def\Pl{\P ^2}

\def\Z{{\Bbb Z}}

\def\W{{\cal W}}

\def\X{{\cal X}}
\def\Y{{\cal Y}}
\def\O{{\cal O}}

\def\N{{\cal N}}

\def\Aut{{\rm Aut}}
\def\Pic{{\rm Pic}}
\def\Spec{{\rm Spec}}
\def\Hom{{\rm Hom}}
 
\def\Im{{\rm Im}}

\def\mult{{\rm mult}}

\def\Yn{Y^{\nu}}
\def\es{e _Y(\Sigma)}

\def\Mg{\overline{M_g}}
\def\Sg{\overline{S_g}}
\def\Hilb{{\rm {Hilb}} ^{p(x)}[\P ^{g-1}]}
\def\PN{P_{N_g}}

\def\HM{\Hom [\Mg ]}
\def\HH{\Hom [H]}
\def\HP{\Hom [\PN ]}

\def\tg3{4{g+1 \choose g-3}}
\def\that{\Theta_{g-3}(Y,W)}
\def\thatt{\Theta_{g-1}(Y,W)}
\def\tdeg{\widehat{\Theta } (Y,W)}
\def\hatt{\widehat{\Theta }}
\def\hats{\widehat{S}}

\cl{\bf{CHARACTERIZING CURVES BY THEIR ODD THETA-CHARACTERISTICS}}

\

\

Lucia Caporaso

Dipartimento di Matematica
Universit\`a degli studi Roma Tre

Largo S. L. Murialdo, 1
I-00146, Roma - ITALY

caporaso@mat.uniroma3.it

\

Edoardo Sernesi 

Dipartimento di Matematica
Universit\`a degli studi Roma Tre

Largo S. L. Murialdo, 1
I-00146, Roma - ITALY

sernesi@mat.uniroma3.it

\

\

\

\cl{\bf  1. Introduction}

\

Consider, for every $g \ge 4$, a non-singular, complex, canonical curve of genus
$g$, that is, a curve $X$  embedded
in $\P^{g-1}$ by its  complete canonical series $|\omega_X|$. Each such a curve
 possesses
$2^{2g}$ line bundles $L$ of degree $g-1$ such that
$L^{\otimes 2}=\omega_X$, the {\it theta-characteristics} of $X$, and precisely 
$N_g := {2^g \choose 2}$ of them are odd (i.e. $h^0(X,L)$ is odd). 
To a non-zero section $\s$ of a theta-characteristic one associates 
a ``half canonical" divisor $D=(\s)$, whose double $2D$ is cut on $X$ by a
hyperplane   $H$ in
$\P^{g-1}$ (whose scheme-theoretic intersection with $X$ is,
of course, everywhere non reduced). For obvious reasons, such a hyperplane $H$
 will be called a {\it theta-hyperplane} of $X$. 

Assume
that $X$ is  general. Then all odd theta-characteristics $L$
satisfy
$h^0(X,L)=1$, while  the even ones have no non-zero sections. Therefore $X$
has exactly $N_g$ theta-hyperplanes; the set of such
hyperplanes will be denoted $\t(X)$, and considered as an element of
$Sym^{N_g}(\P^{g-1})^\ast$.
This said, the main result of this paper is (Theorem \main)

\proclaim Main Theorem. Let $X$ and $X'$ be general canonical curves of
genus
$g
\ge 4$. If $\t(X)=\t(X')$ then $X=X'$.

The case $g=3$, not considered here, has recently been settled in  [CS] (see below).

The above theorem can be put in a different perspective  as follows. Let  $J(X)=Pic^0(X)$  be
the jacobian variety of
$X$, fix a line bundle $L_0$ of degree $g-1$ and let 
$\T=\T_{L_0}=\{\xi\in J(X): h^0(\xi \otimes L_0)>0\}$ 
be the corresponding theta divisor.
If we choose
$L_0$ to be a theta-characteristic we obtain a symmetric theta-divisor. Then the
2-torsion points of $J(X)$ correspond, via multiplication by $L_0$, to the
theta-characteristics on $X$.  If
$X$ is general, $\T$ contains exactly $N_g$ 2-torsion points,
corresponding to the odd theta-characteristics, and they are nonsingular by
Riemann Singularity Theorem. Therefore, by the well known geometric
interpretation of the Gauss map for the theta divisor, we obtain that the
 set of Gauss images of the 2-torsion points of $\T$, call it $\g(J(X),\T _2)$,
coincides with $\t(X)\in Sym^{N_g}(\P^{g-1})^\ast$ 
 (after the choice of a basis of $H^1(\O_X)$).  We  deduce that our theorem
implies the following

\proclaim Theorem.  If $X$ and $X'$ are general abstract curves of
genus $g\ge 4$ such that $\g(J(X),\T_2)$ and $\g(J(X'),\T_2')$ are projectively
equivalent, then $X
\cong X'$.

Observe that this result gives, in particular, a refinement of the classical theorem
of Torelli in the generic case, because $\g(J(X),\T_2)$ only depends on the first
order behaviour at finitely many points of the principally polarized abelian
variety $(J(X),\T)$.  

Notice moreover that the definition of 
$\g(J(X),\T_2)$ can be extended to any principally polarized abelian variety 
$(A,\T)$, provided that the 2-torsion points of the divisor $\T$, which can be
assumed to be symmetric, are all nonsingular. This condition is satisfied on a
dense open subset ${\cal U}$ of the moduli space
${\cal A}_g$ of principally polarized abelian varieties. It is therefore
natural to ask whether  a result  analogous to the above theorem is 
valid for principally polarized abelian  varieties.

This might lead to  an ``odd counterpart'' of the coordinatization of ${\cal A}_g$ 
by  theta constants (see [Mu1] or [Mu3]). We plan to investigate this problem in a
different
 paper.  For a study of similar issues,
 using the classical theta functions, we refer to the work of R.
Salvati Manni in [SM].

Our approach to the Main Theorem is indirect. In fact,
we are  also interested in 
 limits of line bundles and linear series
for families of smooth curves specializing to 
 singular ones. A part of this paper is devoted to certain
aspects of such a  problem. Our results 
are applied to construct
 a degeneration argument that proves our Theorem. 

From this point of view, the first issue is to study
theta-hyperplanes and theta-characteristics for  singular curves,
and to relate the abstract and the projective description.
This is the topic of Section 2. For the abstract question, we use the so-called
moduli space (i.e. stack or scheme)
of spin curves, $\Sg$, constructed by M. Cornalba
in [Co].
$\Sg$ is a geometrically meaningful compactification of the moduli space of
theta-characteristics; it is endowed with a natural, finite morphism
$\pi : \Sg \la \Mg$ 
onto the moduli space of stable curves $\Mg$. The fibers of $\pi$ over
singular curves parametrize their ``generalized theta-characteristics".

On the projective side, the configuration $\t (X)$ of theta-hyperplanes is
naturally defined (in [C]) for certain singular curves in $\P ^{g-1}$;
we thus have a regular morphism
$$
\t :V\la Sym^{N_g}(\P^{g-1})^\ast
$$
where $V$ is a suitable  subscheme of the Hilbert scheme containing canonical
curves in $\P ^{g-1}$. The curves parametrized by $V$ include general smooth curves
and what we call {\it split curves}. A split curve is the
union of two rational normal curves meeting transversely at $g+1$ points. 

In this framework, the crucial result is that split curves 
satisfy   our main
Theorem; that is, a split curve $X_0$ is uniquely determined,
among all curves in $V$, by the configuration
$\t (X_0)$ ([C] Theorem 5, strengthened here by  \AMS).  This is what makes it
possible to prove the Main Theorem by degeneration.

Here is the outline of the argument and a brief description
of the other tools that we
use. Consider a  general 1-parameter 
family $\X \la T$, ($T$ a smooth curve)
of non-singular canonical curves specializing to a split
curve
   $X_0$; denote by $X_t$ the generic fiber. If the Theorem were false,
there would exist a second family   of canonical curves
$
\X'\la  T 
$
with nonsingular generic fiber $X'_t$, and such that
$\t(X_t)=\t(X'_t)$.  
The first difficulty  is to control the special fiber of the second family, 
$X'_0$, over
which there is no a priori information.  In particular, $
\t (X_0')$ might fail to be well defined (if it were defined,  it would obviously
equals
$\t(X_0)$).
To handle the situation,  we study the stable reduction of $\X' \la T$
and 
transfer our projective problem into an  abstract one,
using the existence of $\Sg$ and its properties. 

We proceed to analyze the combinatorial side of the issue;
the goal is to abstractly characterize split curves, by means of their
generalized theta-characteristics; loosely speaking, we are after
an abstract counterpart of the
projective characterization by theta-hyperplanes mentioned above.
This is the topic of  Section 3, which is of independent interest. A purely
combinatorial invariant for an abstract nodal curve (the ``set of exponents") is
defined and it is proved to uniquely determine split curves among all stable
curves (\comb, \split). This set of exponents turns out to encode the numerical
data of the ramification  of the structural morphism $\pi :\Sg \la \Mg$. Whence its
relevance for our central problem.

From the above results we obtain that the stable reduction of $\X' \la T$ has a
 split curve as special fiber. 

Assume now, for simplicity, that $\X' \la T$ is the canonical model of its stable
reduction.
Then $X'_0$ is a split curve and
we can apply the  result mentioned before (\AMS ) to
infer that, since $\t (X_0')=\t(X_0)$, then
$X'_0=X_0$. 

The technical problem, behind the above simplification, is dealt with in Section 4,
by studying the natural action of $PGL(g)$ on  the spaces of configurations of
hyperplanes; we here use methods of  Geometric Invariant Theory.

To conclude that $X_t=X'_t$, 
it suffices to show that the morphism $\t$ is an immersion locally at a  
split curve. This is done in Section 5, 
where we investigate the
tangent space to the locus of deformations of a split curve $X_0$ which remain
tangent to the hyperplanes in $\t(X_0)$. We  apply the theory of
elementary transformations of
vector bundles.

We conclude this introduction with a few words about the case of genus $3$.
In [CS] we proved the Theorem using a far more simple
version of a similar approach. 
In that case, there is no
need to consider  the abstract and the combinatorial side of the
problem. The projective analysis suffices,
using   the techniques of Geometric Invariant Theory.
Just recently, D. Lehavi in [L]  strengthened our result for $g=3$, showing that
every smooth quartic (i.e. not just a general one)
can be recovered by its $28$
bitangents. We thank him for sending us his preprint.

The combinatorial Theorem \comb\  
is a  strengthening of our original statement,
found by  Cinzia Casagrande, to whom we are grateful.

\

\ni
{\bf \NC . Notations and Conventions.}
We work over the field of complex numbers. 
A {\it semistable} curve  in the sense of Deligne and Mumford is a connected,
reduced, projective curve $Y$ having at most nodes as singularities and such that
if
$E\subset Y$ is a smooth rational component, then $E$ meets the union of
the remaining components of
$Y$ in at least $2$ points. 
An $E$ such that $E\cap \overline{Y-E} =2$ is called a {\it destabilizing}
component of $Y$.
If $Y$ is semistable and has no destabilizing component, $Y$ is called a {\it
stable} curve (in the sense of Deligne and Mumford). Stable curves of genus $g$
admit a coarse
 moduli space, denoted by $\Mg$; it is    a projective, integral scheme.

If $Y$ is semistable and no two of its destabilizing components meet, then $Y$ is
called
a {\it quasistable} curve
(a ``decent" curve in the terminology of M. Cornalba, in [Co]) .

An abstract,  stable curve $Y$  of genus $g$ is  a {\it split curve} if
$Y$ is the union of two rational, nonsingular curves meeting transversely 
at $g+1$ points.
A  curve $X\subset \P^{g-1}$ is a  {\it projective split curve} if and
only if
$X$ is the union of two rational normal curves meeting
transversely at $g+1$ points.

Clearly, a projective split curve is the canonical image of an abstract split curve.

We shall consider ``families" 
(of curves, almost always) over a one-dimensional
pointed base
$T$, where $T$  denotes a nonsingular, connected, affine, curve of finite type and
$t_0\in T$ a marked point.
A family ${\cal U}\la T$ is a flat, proper morphism of schemes; the fiber over the
marked point $t_0$ will be denoted  by $U_0$ or simply by $U$, and called the
``special" fiber. The fiber over $t\neq t_0$ will be denoted bt
$U_t$ and called the ``generic" fiber. 

In certain contexts, a family
${\cal U}\la T$
as above will be called a {\it one parameter deformation} of its central fiber
$U_0$.

We shall denote $\PN :=Sym^{N_g}(\P^{g-1})^\ast$.

For $S$  a topological space, $\g _c(S)$ is the number of its
connected components.

\

\ni
{\bf \PO.} 
Let $L$ and $M$ be two finite sets
of integers.
We say that $L$ {\it dominates} $M$ (in symbols $L \geq M$) iff there
exists a surjective map 
$
\a : L\la M
$
such that for every $l \in L$ we have $\a(l) \geq l$.
It is  very easy to see  that,
if $L\geq M$ and $M\geq L$, then $L=M$. 
In particular, the above definition
 gives a partial ordering on the set of all
finite sets of integers.

Let $S$ be a purely dimensional scheme (that is, every irreducible component
of
$S_{red}$ has the same dimension). We associate to $S$ its {\it multiplicity
set}
$L(S)$  as follows
$$
L(S) :=\{ n : \exists \  {\rm{some \  irreducible \  component}}\  Z\   {\rm of}
\  S\   {\rm such \  that }\  \mult _Z S = n \}
$$

\

\cl{\bf 2. The moduli theoretic framework}

\

\ni
{\bf \BS . The basic projective set up.} Fix $g\geq 3$ and consider the  
Hilbert scheme $\Hilb$ of curves in $\P ^{g-1}$ having Hilbert polynomial 
$p(x)=(2g-2)x-g+1$.
In it we find the locus  of connected curves $X$ 
of degree $2g-2$ and arithmetic genus $g$, satisfying the
following three conditions:

\

\ni
{\bf 1.} $X$ is reduced and has at most 
nodes as singularities.

\ni
{\bf 2.} $X$ is embedded in $\P^{g-1}$ 
by the complete linear series $|\omega_X|$ (where $\omega_X$ denotes the dualizing
line bundle of $X$).

\ni
{\bf 3.} No irreducible component of $X$ is contained in a hyperplane.

\

\ni
A curve satisfying 1. and 2. above is, clearly, stable in the sense of Deligne
and Mumford; we shall call it
 a {\it canonical} curve throughout the paper. 

\

\ni {\bf Definitions.} Let $X\subset \P^{g-1}$ be a  curve satisfying 1, 2 and 3 
above.  A hyperplane $H\subset \P^{g-1}$ will be called a {\it
theta-hyperplane} of
$X$ if $H\cap X$ is everywhere non reduced and  supported at $g-1$ distinct
points. Let $i$ be an integer with $0 \le i \le g-1$. We shall say that the  
 theta-hyperplane $H$ is of type $i$  if $H$ contains exactly $i$ singular
points of $X$. 
Such a curve $X$ will be called {\it theta-generic} if it has finitely many
theta-hyperplanes.

\

\ni
As noticed in the introduction, a general nonsingular canonical curve  
has  $N_g = 2^{g-1}(2^g-1)$ distinct theta-hyperplanes,
corresponding bijectively to its odd theta characteristics.
 Therefore it is theta-generic. 

We shall say that an abstract stable curve  is theta-generic if its dualizing line
bundle is very ample and if its canonical model  is theta-generic.
Notice that theta-generic curves can be of only two topological types:
either they are irreducible, or they are the union of two rational components
meeting at $g+1$ distinct points, i.e. they are split curves.

We shall denote by $V \subset \Hilb$ the open subset parametrizing 
theta-generic  curves and by $V_0\subset V$ the open subset 
parametrizing nonsingular theta-generic  curves. 
As in [C] and [CS] we  define a morphism
$
\t: V_0 \to Sym^{N_g}(\P^{g-1})^\ast
$
by  $\t(X) = \{H_1, \ldots, H_{N_g}\}$, where $H_1, \ldots, H_{N_g}$
are the (distinct) theta-hyperplanes of $X$. 
Then, arguing as in [CS], Lemma 2.3.1, it is easy to see that   
$\t$ can be extended to a morphism (called again $\t$)
$$
\t: V \to Sym^{N_g}(\P^{g-1})^\ast
$$

It will often be convenient  to view 
$\t(X) $ as a (not necessarily reduced)
hypersurface
of degree
$N_g$ in
$\P^{g-1}$, all of whose irreducible components are (possibly multiple)
hyperplanes; such a hypersurface will be also denoted $\t(X)$, by abuse of 
notation, and called the {\it theta-hypersurface}
of $X$.
From [C] and [CS] we get that $\t (X)$ is reduced if and only if $X$ is smooth;
if $X$ is singular, all hyperplanes of type $i$ appear in $\t (X)$ with
multiplicity $2^i$.

For a given $X$ in $V$ we will denote by $t_i(X)$ the number of distinct
theta-hyperplanes of type $i$. The numbers $t_i(X)$ only depend on the
number of nodes and of irreducible components of $X$, and have been computed in [C].

Let 
$$
\matrix{\X & \hookrightarrow {} & T\times \P^{g-1}\cr
\   \lmapdown{ } & &\cr
T& &}
$$
be a family of    curves whose  Hilbert polynomial is $p(x)$;
the natural morphism to the Hilbert scheme is denoted by 
  $$\psi _{\X}:T\la \Hilb .$$  Assume that the generic fiber is in $V$,
 that is $\psi _{\X}(T-\{t_0\})\subset V$.
Since $T$ is a nonsingular curve and $\PN:=Sym^{N_g}(\P^{g-1})^\ast$ is projective
the composition  $\t\circ\psi _{\X}$ extends to the whole of $T$; it will be
denoted by
$$
\t _{\X}:T \la \PN .
$$ 
If $\psi _{\X}(T-\{t_0\})\subset V_0$ we consider the curve
$$
J_{\X}\subset T\times (\P ^{g-1})^*
$$
defined as the closure of the incidence correspondence
$$\{ ( t,H):\  H\subset \t (X_t), \  t\neq t_0\}\subset T\times (\P ^{g-1})^*$$
Away from $t_0$, $J_{\X}\la T$ is an unramified  covering of degree $N_g$.

In the next section, we describe the abstract counterpart of the above set-up.

\

\ni
{\bf \CS . Cornalba's moduli space of spin curves.} In [Co] M. Cornalba constructed
a geometrically meaningful compactification
$\Sg$, over $\Mg$, of the moduli space of theta-characteristics of smooth curves
of genus
$g$. 
We need to recall some basic features of $\Sg$, which is called the
{\it moduli space of stable spin curves}.

Both the scheme description and the stack description are available, we shall
confine ourselves to the schematic definition, which suffices for our purposes.

$\Sg$ is a normal, projective scheme  which admits a
proper morphism $\pi$ onto $\Mg$
$$
\pi :\Sg\la \Mg.
$$

As expected, the degree of $\pi$ is $2^{2g}$, and
$\Sg$ is a disjoint union of two
irreducible components,
$\Sg^+$ and
$\Sg  ^-$,  corresponding, respectively, to even and odd theta-characteristics.
We are mostly interested in odd theta characteristics in this paper.
The degree of the restriction of $\pi$ to $\Sg ^-$ is $N_g$.

The points of $\Sg$ are described as line
bundles on certain quasistable 
curves of arithmetic genus $g$.
Some more notation: let $Y$ be  a stable curve  and
let $\Sigma$ be a set of nodes of $Y$, denote  by $\Yn_{\Sigma}$ the
normalization of $Y$ at all nodes in $\Sigma$ and by  $Y_{\Sigma}$ the
quasistable curve   obtained by ``blowing-up" $Y$ at all nodes in $\Sigma$.
Thus  $\Yn_{\Sigma}$ is  the closure of what remains of
$Y_{\Sigma}$  after removing all
of its destabilizing components.

A point in $\Sg$  is (the isomorphism class of) a {\it spin curve} $\xi$,
that is the following set of data.
First, a subset $\Sigma _\xi = \Sigma$ of nodes  of $Y$ and the quasistable
curve
$Y_{\Sigma}$. Second, a line bundle $L_\xi =L$ on $Y_{\Sigma}$ and a
homomorphism
$\a:L^2\la \omega _{Y_{\Sigma}}$ such that
$L$ has degree $1$ on all destabilizing components, and the restriction of 
$\a$
 to 
$\Yn_{\Sigma}$
is an isomorphism (see [Co] for details).
 We shall say that
$\xi$  is {\it{supported on}} $\Sigma$.
As it will be made clear later in the paper (Section 3),
not all subsets of nodes of $Y$  are
the support of a spin curve: those that are will be called ``admissible"
and will be studied extensively later.

The scheme structure on $\Sg$ is obtained  in [Co] after
constructing the ``universal deformation space" of a spin curve $\xi$.
First (in Section 2) he describes its group of automorphisms $\Aut (\xi )$ and shows
that there is an exact sequence of finite groups
$$
1\la \Aut _0(\xi )\la \Aut (\xi ) \la \Aut (Y)
$$
where $\Aut _0(\xi )$ is the so-called group of ``inessential" automorphisms,
(acting trivially on $Y$).
Call $B_{\xi}$ the base of the universal deformation space of $\xi$ and
$B_Y$ the base of the universal deformation space of $Y$; $B_\xi$ and $B_Y$
will be viewed, as in [Co], as analytic spaces, that is, as $3g-3$-dimensional
discs. There is a commutative diagram of morphisms
$$
\matrix{B_{\xi} & \mapright {} & B_{\xi}/\Aut _0(\xi) 
&\mapright {}& B_{\xi}/\Aut (\xi)
\  &\hookrightarrow {}&\Sg\cr
\lmapdown{\d } & & \lmapdown{\pi _{\xi} } & & & &\mapdown {\pi }\cr
 B_Y &= &   B_Y \  \  &\mapright {\rho _Y}&  B_Y/\Aut (Y)
\  &\hookrightarrow {}&\Mg }
$$
with $B_{\xi}/\Aut _0(\xi)  \subset B_Y\times _{\Mg}\Sg $; everything is shown to
satisfy the necessary compatibility conditions. We recall how the covering $\d$ is
defined ([Co] section 5): choose coordinates
$u_1,....,u_{3g-3}$ for ${ B_Y}$ so that the first $\#\Sigma$
correspond to the loci where the $i$-th node of $\Sigma$ is preserved.
Choose cordinates $t_1,....,t_{3g-3}$ on $B_{\xi}$ and define $\d$ as the base
change 
$u_i=\d ^*(t_i^2)$ for $i\leq \#\Sigma$ and $u_i=\d ^*(t_i)$ for $i>\#\Sigma$.

Let $\Y \la T$ be a family of generically stable curves.
Since $T$ is nonsingular and $\Mg$ is projective, its  moduli map 
$T - \{t_0\} \to \Mg$ extends 
to $T$ and is denoted by
$$\phi _{\Y} :T \la \Mg$$
The pull-back to $T$ of $\Sg ^- $ is a curve over $T$
denoted by
$$
S_{\Y}:=\phi _{\Y}^*\Sg^-.
$$

\proclaim Proposition \exp  .
Let $\Y \la T$ be a  general one-parameter deformation of the stable curve
$Y$. 
\item a ) The  curve $S_{\Y }$ is smooth.
\item b) Let $\xi$ be a spin
curve on $Y$ supported on the set
$\Sigma$. 
Then the index of ramification of 
the finite covering \  $S_{\Y }\la T$ at the point corresponding to
$\xi$ is
$2^{e_Y(\Sigma)}-1$
where $e_Y(\Sigma)= \#\Sigma - \g _c(\Yn_{\Sigma}) +1 $.

\rk The same result holds for even spin curves, that is, if we replace
$S_{\Y}$ by $\phi _{\Y
}^*\Sg ^+$. The proof is the same.

\

\pf 
It suffices to show that  $S_{\Y}$ is non singular at
the points lying over $t_0$.
The moduli morphism $\phi _{\Y}$ factors locally
through the  map $\phi _{Y}:{\cal O}_{T,t_0}\la B_Y$;
let ${\overline T}:=\Im\phi _{Y}$.
We shall prove, more precisely,
 that if  ${\overline T}$ 
is  transverse in $B_Y$  to the loci where the
nodes of $Y$ are preserved, then $S_{\Y}$ is smooth.

It  suffices to prove
that
$\phi _{Y}^*\circ\rho _Y^*\Sg ^-$ is smooth. By the above diagram, this is
equivalent to show that the curve $\pi _{\xi} ^{-1} {\overline T}$ is smooth.
Such a curve is the quotient via a finite group of the curve $\d ^{-1}{\overline
T}$; this last curve is smooth,
because, by assumption, it is transverse to the branch locus of $\d$  (look at
the explicit description of
$\d$ given above). Hence its quotient by any finite group is smooth; this
proves a).

For b) as well, 
 most of the argument  
is already in Cornalba's paper (Section 5), together with various explicit
examples (from (5.3) to (5.7) ). 

It suffices to look at how the covering $\pi_{\xi}^{-1}({\overline T})\la
{\overline T}$ ramifies.
The
natural map
$\d :B_{\xi}\la  B$  is a finite covering of degree $2^{\#\Sigma }$
totally ramified over the origin (as explained above);
$\pi_{\xi}^{-1}({\overline T})$ is the quotient of $\d ^{-1} ({\overline T})$
via $\Aut _0 (\xi )$.
The structure of 
$\Aut _0(\xi)$  is described in [Co], Lemma 2.2 
in terms of the  graph associated to $
\xi$
(having verteces the connected components of $\Yn _{\Sigma}$ and edges the
destabilizing components of $Y_{\Sigma}$); his result can be re-stated by saying
that 
$\Aut _0(\xi)$ is a vector space over $\Z /2\Z$ 
of dimension equal to the number of connected components of $\Yn _{\Sigma}$ minus $1$:
$$
\dim _{\Z /2\Z}\Aut _0(\xi) = \g _c(\Yn _{\Sigma }) -1.
$$
 Hence $\pi_{\xi}^{-1}({\overline T})\la
{\overline T}$  ramifies  at the point corresponding to $\xi $ with index
$2^{\#\Sigma -\g_c (\Yn _{\Sigma }) -1}+1$.
\qed

\

\ni
{\bf \rig .} We now introduce the scheme $S_Y$ parametrizing odd spin curves having $Y$ as
stable model:
$S_Y$ is the (scheme-theoretic) special fiber of $\rho _Y^*\Sg ^-=B_Y\times
_{\Mg} \Sg ^- \la B_Y$.

Thus $S_Y$ is a zero-dimensional scheme of length $N_g$. If $\Aut(Y)$ is trivial,
then $S_Y$ is the fiber of $\Sg ^-$ over the point $Y\in \Mg$.
The following statement is an obvious consequence of \exp.
 \proclaim Corollary \lsy. Let $Y$ be a stable curve. $L(S_Y)=\{ 2^n,\  {\rm
s.t.}\ 
\exists
\xi\in S_Y : e_Y(\Sigma _{\xi} )=n\}$.

  The number $\es$, defined in the statement of \exp, will
be called the {\it exponent} of $\xi$ and will play a crucial role in the
sequel. Notice that, denoting by $Z=\Yn _{\Sigma}$ and by $g_Z$ its arithmetic
genus, we have that 
$\es = g - g_Z$. The following simple observation will be used later. 

\proclaim Claim. Assume that $\xi \in S_Y$ has exponent at least $g-3$;
then  $h^0(Z,L_\xi\otimes {\cal O}_Z ) =1$.

In fact $\xi$ is an
odd spin curve, by definition;
thus $L_\xi\otimes {\cal O}_Z$ is an odd (hence effective) square root of  
$\omega _Z$,
 that is, $h^0(Z,L_\xi\otimes {\cal O}_Z )$ is odd (see [Co] section 6).
Now $Z$ is a curve of genus at most $3$,  hence 
of course $h^0(Z,L_\xi\otimes {\cal O}_Z)\leq 2$. 

\

Denote by $S_{\xi}
\subset S_Y$  the (unique) connected component of
$S_Y$ supported on $\xi$. Denote by $\hats _Y$ the union of all components of $S_Y$
corresponding to spin curves whose multiplicity is at least $2^{g-3}$, that is
$$
\hats _Y:=\bigcup_{\xi : e_Y(\Sigma _{\xi})\geq g-3} S_{\xi}
$$
Of course, $\hats _Y$ can be empty, for example, it will be empty if $Y$ is of
compact type, in which case $S_Y$ is reduced (see [Co] and also  \ctp). 

\

\ni
{\bf \AP . Relating the abstract and the projective pictures.}
The set up in this section is the following. $Y$ is a stable curve and $\Y \la
T$ a one-parameter deformation of $Y$ with smooth and theta-generic general
fiber.  A {\it  generically canonical, birational
model} of $\Y \la T$ is the following: a projective family  of curves:
$
\P ^{g-1} \times T \supset \W \la T
$
and a birational $T$-map $\rho:\Y --> \W$ which is an isomorphism away from
$t_0$ and such that for $t\neq t_0$, the induced map $Y_t\la W_t$ is the
canonical morphism.

Unless otherwise specified, no assumption will be made on the central fiber $W$
of $\W$, or on
the restriction of
$\rho$ on the central fibers (notice that the locus of canonical curves is not
closed in $\Hilb$).

To such a picture we add two, differently defined, families of odd
theta-characteristics: to the abstract family $\Y \la T$ we associate $S_{\Y}\la
T$ (defined in \CS); to the projective family $\W \la T$ we associate $J_{\W}\la T$
(defined in \BS).
To make our notation more precise, consider
$\theta _{\W}:T\la \PN$, defined in \BS. For $t\neq t_0$, $\theta
_{\W}(t)$ represents the configuration of the $N_g$ theta-hyperplanes of $W_t$.
Since $W_t\in V$, this configuration does not depend on $\W$ and we have
$\theta
_{\W}(t) =\theta (W_t)$. For $t_0$ the configuration $\theta _{\W}(t_0)$
corresponds to a hypersurface of degree $N_g$ in $\P ^{g-1}$ for which we shall
use the notation
$$
\theta ^{\W}(W):=\theta _{\W}(t_0)
$$
(as usual, the above symbols will be abused to denote both   points in $\PN$
and  hypersurfaces in $\P ^{g-1}$). If $W\in V$ we have, of course, $\theta
^{\W}(W)=\t(W)$ defined before.
Consider now $J_{\W} \subset T\times (\P
^{g-1})^*$. 
If $t\neq
t_0$, the fiber of $J_{\W}$ over $t$ only depends on $W_t$
  and will be denoted   by  $J_{W_t}$.
The special fiber $(J_{\W})_{t_0}$   depends, a priori, on $\W$. We shall denote
$$
J_W^{\W}:=(J_{\W})_{t_0}
$$
and, if $W\in V$ we shall simply write $
J_W:=(J_{\W})_{t_0}
$.

Let us consider  the multiplicity
sets;  if $t\neq t_0$ we have, obviously,
$
L(\theta (W_t)) = L(J_{W_t})=\{ 1\}
$. 
For the special fibers
we have, by construction,
$
L(\theta ^{\W}(W)) = L(J_W^{\W}).
$

\proclaim Lemma \cor .
Let $\Y \la T$ be a  general one-parameter deformation of the stable curve
$Y$;  let $ \W \la T$ be a generically
canonical, birational image of
$\Y
\la T$.
(a) There  exists a natural, surjective, birational $T$-morphism
$$
\mu : S_{\Y} \la J_{\W}.
$$
(b) Assume furthermore that $K_Y$ is very ample and that the central fiber $W$
of $\W \la T$ is a canonical image of $Y$, with $W\in V$.
Then the morphism $\mu$ above induces an isomorphism
on the central fibers.
In particular, $
L(S_Y)=L(J_W).
$

\pf By assumption, the generic fiber $Y_t$ of
$\cal Y$ is  smooth and non-hyperelliptic  and the generic fibers of  $\W\la T$
are in
$V$.
There
is a natural
$T$-birational map
$$
\mu  :S_{\Y} ---> J _{\W}
$$
which is an isomorphism away from $t_0$.
It associates to an odd theta-characteristic on
$Y_t$ the corresponding theta-hyperplane of its canonical model $W_t$.
By part (a) of  Proposition \exp \  the morphism $\mu $ extends to the whole of
$S_{\Y}$. 
This proves (a).

For (b), notice that by (a) it is enough to show that $\mu$ induces a set
theoretic bijection on the central fibers,
that is the set theoretic map
$$
\mu _0 :(S_Y)_{red} \la (J_W)_{red}
$$
is a bijection.
The proof is a bit long, but very simple and standard.

Throughout the rest of the argument, we shall identify $Y$ with $W$ 
by the canonical isomorphism given in the statement.
Let $\xi\in S_Y$ be a spin curve supported on the set of nodes
$\Sigma _{\xi}
=
\Sigma$. To
$\xi$ there corresponds a line bundle $L_{\xi}$ on the quasistable curve
$Y_{\Sigma}$ such that, if $Z:=\Yn _{\Sigma}\subset Y_{\Sigma}$, we have
$(L_{\xi}\otimes {\cal O}_Z)^2\cong \omega _Z$ and such that,
for every destabilizing component $E$ of $Y_{\Sigma}$, the restriction of
$L_{\xi}$ to $E$ is ${\cal O}_E(1)$.
Since $W\in V$, there exists a unique effective 
divisor
$D_{\xi}\in \Pic Z$ such that
$$
{\cal O}_Z(D_{\xi}) = L_{\xi}\otimes {\cal O}_Z.
$$
Moreover, $D_{\xi}$ is reduced and $suppD_{\xi}\subset (Y_{\Sigma})_{smooth}$.

Consider the canonical morphism $\sigma : Y_{\Sigma}\la W\subset \P ^{g-1}$
(which of course contracts all the destabilizing components); locally at every
point of $D_{\xi}$, $\sigma $ is an isomorphism.
The set of points $\Sigma \cup \sigma (D_{\xi})\subset W$ is a set of $g-1$
points in general position (recall that $W$ is theta-generic), which therefore
spans a unique hyperplane $H_{\xi}$ containing $\Sigma$.
It is clear that $H_{\xi}$ is tangent to $W$ at every point of $\sigma
(D_{\xi})$.

In conclusion, we have explicitely described $\mu _0(\xi)$; this has the
advantage of proving that $\mu _0$ is injective. In fact
$\mu _0(\xi )=\mu _0(\xi ') \Leftrightarrow H_{\xi} = H_{\xi '}$,
this implies that $\Sigma _{\xi } = \Sigma _{\xi '}$ and that $D_{\xi } =
D_{\xi '}$. This is enough to ensure that $\xi = \xi '$ (see [Co]).

It remains to show that $\mu _0$ is surjective.
Let $H$ be a theta-hyperplane of $W$, so that $(t_0,H)\in J_W$, and let $\Sigma
= H\cap W_{sing}$. Then $H$ is tangent to $W$ at $g-1-\#\Sigma$ smooth points of
$W$, let $D$ be the divisor defined by $H\cap W_{smooth} = 2D$.
Denote by $\sigma :Y_{\Sigma }\la W$  the canonical morphism, let
$Z:=\Yn _{\Sigma}\subset Y_{\Sigma}$ and 
denote by $\sigma _Z :=\sigma _{|Z}$ the restriction to $Z$.
Let 
$\Delta :=\sigma ^*_Z(\Sigma )$ so that $\deg \Delta = 2 \#\Sigma$ and 
$$
\sigma _Z^*H=2\sigma _Z^*(D )+\Delta .
$$
On the other hand, $W$ is a canonical image of $Y_{\Sigma}$ hence
$$
\sigma _Z^*H=\omega _{Y_{\Sigma}}\otimes {\cal O}_Z=\omega _Z\otimes {\cal
O}_Z(\Delta ).
$$
Hence $\omega _Z =2\sigma _Z^*(D )$.

Finally, by gluing $\sigma _Z^*(D )$ to ${\cal O}_E(1)$ on every
destabilizing component $E$ of $Y_{\Sigma}$, we obtain a line bundle $L$
which 
(regardless of the gluing data, which are proven to be irrelevant in [Co])
corresponds to a spin curve $\xi$ on $Y$ such that $\mu _0(\xi) = H$.
Proving that $\mu _0$ is surjective.
\qed

Let $H\subset \Hilb$ be the projective scheme defined as the closure  of $V$.

Let $P$ be a projective scheme and consider the space of morphisms from $T$ to $P$
$$
\Hom [P]:=\{\tau :T \la P\}.
$$
Since $T$ is a smooth curve,  any rational
map from $T$ to  $P$ extends to a (uniquely defined) regular
morphism 
${\overline T}\la P$, where ${\overline T}$ is the  smooth compactification of $T$. 
Thus $\Hom [P]$ has a natural scheme structure.
We
can apply this to the projective varieties $\Mg$, $H$ and
$\PN$.
We shall thus consider the scheme of maps from $T$ to $\Mg$:
$
\HM 
$; 
  the scheme of maps from $T$ to $H$:
$
\HH 
$,
and the scheme
$
\HP 
$ as above.

There is a rational  map $\phi  :H --> \Mg$ (regular at least on $V$),
which is the moduli map associated to the universal family over $H$
Recall also that we have a rational map $\theta : H --> \PN$ which is regular
on $V$.
 Therefore we have two morphisms of schemes:
$$
\phi _* :\HH \la \HM
$$
given by composing with $\phi$ (that is, for $\psi \in \HH$, we define 
$\phi _* (\psi ):=\phi \circ \psi$);
$$
\theta _*  :\HH \la \HP
$$
given by composing with $\theta$ (as above).
We have

\proclaim Lemma \num .
Let $\Y \la T$ be a one-parameter deformation of $Y$, with smooth
and theta-generic general fiber.
Let $\W \la T$ be a generically canonical, birational image of $\Y$. Then
$$
L(S_{Y})\geq L(J_W^{\W}).
$$

\pf
First we consider the case in which $\Y \la T$ is a general one parameter
deformation of $Y$.
By  \exp, we have that $S_{\Y}$ is a smooth curve  and,
by  \cor (a),  there
is a surjective $T$-morphism 
$$
S_{\Y}  \la J_{\W}.
$$
By comparing the fibers over $t_0$ we obtain
$
L(S_{Y} )\geq L(J_W^{\W})
$
which is what we wanted.

Consider now a one-parameter deformation $\Y\la T$ of $Y$
and a generically canonical model $\W\subset T \times \P ^{g-1}$
 as  stated. 
As usual, denote by  $\psi _{\W}\in \HH$ and $\phi _{\Y}\in \HM$ the
corresponding moduli morphisms; of course, $\phi _{\Y}=\phi\circ \psi
_{\W}$.
Pick now a local  curve $U^M \subset \HM$ whose special point is $\phi
_{\Y}$. More precisely, let $R$ be a discrete valuation ring and let
$$
U^M:=\Im \{ \phi _R:\Spec R \la \HM\}
$$
be such that the image of the generic point,
$\phi _R(\eta)=\phi _\eta$, is a general
one-parameter deformation of
$Y$, and the image of the special point is the given $\phi _R(s)=\phi
_{\Y}$.

Let $U^H\subset \HH$ be a lifting of $U^M$ passing through $\psi _{\W}$; that
is, 
$$
U^H:=\Im \{ \psi _R:\Spec R \la \HH\}
$$
such that 
the image of the generic point,
$\psi _R(\eta)=\psi _\eta$ is such that $\phi\circ \psi _{\eta} = \phi _{\eta}$
while the image of the special point is the original $\psi _{\W}$.
Denote by $\W ^{\eta}\subset T\times \P^{g-1}$ the family of curves
corresponding to
$\psi _{\eta}$ (that is, $W^{\eta}$ is the pull-back, via $\psi _{\eta}$, of the
universal family over
$H$ and $\psi _{\eta} =
\psi _{\W ^{\eta}}$). 

Lastly, let $U ^P\subset \HP$ be defined as $U^P:=\theta _*(U^H)$.
Thus 
$$
U^P= \Im \{\theta _R : \Spec R \la \HP \}
$$
such that, denoting, as usual, $\theta _{\eta}$ and $\theta _s$ the images of
the generic and special point (respectively) of $\Spec R$, we have
$
\theta_{\eta} = \theta \circ \psi _{\W^{\eta}} = \theta _{\W ^{\eta}}
$
and 
$
\theta _s = \theta _{\W}.
$
Now we shall construct a universal family of theta-hypersurfaces over $U=\Spec
R$ using the morphism defined above
$$
\t _R :U \la \HP \hookrightarrow Hilb (T\times \PN )
$$ 
(where the  immersion on the right is the structural one, used to
give a scheme structure to $\HP$).
To do that, we shall proceed in a standard
way, using the various universal families over the Hilbert schemes in our
set-up.  

There is a universal family  
${\cal F}\la Hilb (T\times \PN)$
 which we can pull back to $U$ via $\t _R$. 
Denote by  ${\cal G} := {\cal F} \times _{Hilb (T\times \PN )}U\la U$ this
pull-back. All of its fibers are thus isomorphic to $T$. Since $U=\Spec R$, we get
that this fibration is in fact trivial, so that
$
{\cal G}\cong T\times U;
$
let us fix one of such isomorphisms from now on.

Consider also the universal family  of 
(reducible) hypersurfaces over
$\PN$, denoted by  $\PN\times \P^{g-1}\supset {\cal I}\la \PN$.
We  have a diagram
$$
\matrix{{\cal G}&\la &{\cal F}  &\hookrightarrow  Hilb (T\times \PN )\times
T\times \PN
\la \PN\cr
\lmapdown {}&&\lmapdown {}&\cr
U &\mapright{} &  Hilb (T\times \PN ) &}
$$
so that we can pull back ${\cal I}\la \PN$
to ${\cal G}$. Denote $\widetilde{\Theta}:={\cal G}\times _{\PN}{\cal I}$ such
a pull back. Finally
$$
\widetilde{\Theta}  \la {\cal G}\cong T\times U\la
U 
$$
is the  family of theta-hypersurfaces that we wanted.
Let us now consider the restriction ${\widetilde \T}^0$ of ${\widetilde \T}$
to $\{t_0\}\times U$, i.e. 
$
{\widetilde \T}^0 \la \{t_0\}\times U \cong \Spec R
$.
This is a 
family of schemes of pure dimension $g-2$, for which
we get
$
L(({\widetilde \T}^0 )_{\eta}) \geq L(({\widetilde \T}^0 )_s).
$
By construction,  the special fiber is the original
$\t^{\W}(W) =({\widetilde \T}^0 )_s  $.
For the generic fiber, of course,
$({\widetilde \T}^0 )_{\eta}=\theta ^{\W ^{\eta}}(W^{\eta})$
(where, as usual, $W^{\eta}$ denotes the central fiber of $\W^{\eta}$). 
Thus the above inequality translates into
$
 L(\theta ^{\W ^{\eta}}(W^{\eta}))\geq L(\t^{\W}(W)).
$
Recall
now that
$\W ^{\eta}$ is a birational, canonical image of a general one-parameter
deformation
 of $Y$. Therefore we can apply to $\W ^{\eta}$ the result of the  first
part of the proof, that is:
$
L(S_{Y})\geq L(J^{\W^{\eta}}_{W^{\eta}});
$
since
$
L(J^{\W^{\eta}}_{W^{\eta}}) = L(\theta ^{\W ^{\eta}}(W^{\eta}))$,
we conclude that
$
L(S_{Y})\geq L(\t^{\W}(W)) = L(J^{\W}_W)
$.
\qed

\

\cl{\bf{3. Combinatorics of stable curves}}

\

This  section is essentially  independent of the rest of
the paper; its goal is \split, 
which will play a crucial role
later on.
We have seen that, for a stable curve $Y$, the multiplicity set $L(S_Y)$ only
depends on the combinatorial data of $Y$ (\exp \  and  \lsy).
Here we shall make this more precise by re-defining it in purely combinatorial
terms and thus making its computation rather simple.
The key point is that $L(S_Y)$ is a fine invariant for certain types of curves:
 namely, 
it completely determines  the
topological and combinatorial structure of $Y$
(among all stable curves)   if $Y$ is of compact type
(this is easy: \ctp) and if $Y$ is split (this is more interesting: \comb,
\split).

\

\ni
{\bf \CI. } Let $Y$ be a  nodal curve.
$\Gamma _Y$ is the dual graph of $Y$,
that is, the verteces of $\Gamma _Y$ correspond to the irreducible components
of $Y$, the edges joining two verteces correspond to the nodes contained in the
two corresponding components.

Recall that $Y_{sing}$ is the set of all nodes of $Y$.
If $\Sigma\subset Y_{sing}$ the curve $Y^{\nu}_{\Sigma}$ is the
normalization of $Y$ at the nodes contained in $\Sigma$.
Assume that $Y$ has arithmetic genus $g$,  $\d$
nodes and
$\g$ irreducible components; let 
$Y=\bigcup _{i=1}^\g C_i$
 be the
decomposition of
$Y$ into irreducible components, and let $g_i$ be the geometric genus of $C_i$.
Recall the genus formula: 
$
g=\sum g_i +\d -\g +\g _c(Y).
$

 Let $N$ be a node of $Y$. We say that $N$
is {\it internal} if
$N$ is contained in a unique irreducible component of $Y$, otherwise we say that
$N$ is {\it external}. A (external)  node $N$ is called {\it separating}  if
$\g _c(Y) = \g _c(Y- N) -1$. If all nodes of $Y$ are separating, $Y$ is called
``of compact type".

\
 
\df {\it Exponent of a set of nodes.}
Let $\Sigma$ be a set of nodes of $Y$. The
{\it{exponent}} of
$\Sigma$ is  the number $e_Y(\Sigma)$ below
$$
\es = p_a(Y)-p_a (\Yn_{\Sigma})
$$
As we already know (see \rig), 
an equivalent definition is
$
e_Y(\Sigma):= \#\Sigma - \g _c(\Yn_{\Sigma})
+\g _c(Y).
$

The reason for the name ``exponent"
comes from Proposition \exp  . 
Keeping the  same notation, we list a couple of straightforward consequences of the
definition.

\proclaim Property A.
If $\Sigma \subset \Sigma '$ then $e_Y(\Sigma) \leq  e_Y(\Sigma ')$.
Moreover, for every $\Sigma$,
$$
0=e_Y(\emptyset )\leq e_Y(\Sigma) \leq e_Y(Y_{sing})=g-\sum _{i=1}^{\g}g_i
$$

\proclaim Property B.
If $N$ is a node of $Y$ not contained in $\Sigma$,
then

\centerline{$
e_Y(\Sigma \cup \{N\})=\cases{e_Y(\Sigma ) +1 & {\it if \ $N$ is not separating
for $\Yn _{\Sigma}$;}\cr
e_Y(\Sigma ) & {\it if \  $N$ is  separating for $\Yn _{\Sigma}$.}\cr}
$}

\

\ni
We shall  need to compute the exponent of a distinguished type of sets of nodes:

\

\df {\it Admissible sets of nodes.}
A set $\Sigma \subset Y_{sing}$ is said to be {\it admissible} 
if for every subcurve
$W$ of
$Y$,
the number of nodes in the intersection $W \cap
W^c$ which are not in
$\Sigma$ is even. 
Equivalently, $\Sigma$ is admissible if for every $i=1,....,\g$ the number
$\#[(C_i\cap C_i^c)\cap(Y_{sing} -\Sigma)]$ is even.

\

The equivalence of the two definitions above is an elementary verification.
We denote by $A_Y$ the set of all admissible sets of nodes of $Y$:
$$
A_Y:=\{ \Sigma \subseteq Y_{sing}\  \rm{s.t. \  } \Sigma \  \rm{is \ 
admissible}\}
$$

The previous definition is implicit in Cornalba's paper;
it is
motivated by the fact
that  
$\Sigma$ is admissible if and only if  the  dualizing sheaf $\omega _{\Yn
_{\Sigma }}$ has even degree on every subcurve. 
It  is important to stress that a set $\Sigma$ is the support of a spin
curve on
$Y$ if and only if $\Sigma \in A_Y$ (see [Co]).

Some simple examples; the set of all
nodes
$\Sigma = Y_{sing}$ is certainly admissible.  Similarly the set of all external nodes
of $Y$ is admissible, in fact the above definition ignores the internal nodes. 
The empty set is not always admissible, consider for example a singular curve of compact type
$Y=C_1 \cup C_2$ with $\#C_1\cap C_2 =1$.

Here are a few 
immediate properties:

\proclaim Property C.
If $\Sigma \in A_Y$, then 
$\Sigma$ contains all the separating nodes of $Y$.

\proclaim Property D.
Let $\Sigma \subset Y_{sing}$, denote $Z=\Yn _{\Sigma}$ and let $\Sigma
_1\subset Z_{sing}$. Then $\Sigma _1\in A_Z $ if and only if $ \Sigma _1\cup
\Sigma\in A_Y $. In particular $\emptyset \in A_Z $ if and only if $
\Sigma\in A_Y $.

\

\ni
{\bf \cyc. Cyclic sets of nodes.}
Given our nodal curve $Y$, consider its dual graph $\Gamma _Y$.
To any subset, $\Sigma$, of nodes of $Y$,
 we can associate a graph $\Gamma
_{\Sigma}$, which is the subgraph of $\Gamma _Y$ generated by all edges 
representing the nodes contained in $\Sigma$. We shall say that a non-empty $\Sigma$ is a {\it
cyclic} set of nodes if its  graph $\Gamma
_{\Sigma}$ is a closed polygon. 

An equivalent, graph-free, definition is the following.
Up to reordering the irreducible components of $Y$, we can write a cyclic set
of nodes as
$$
\Sigma =\{N_{1,2},N_{2,3},...,N_{h-1,h},N_{h,1}\}
$$
with $1\leq h\leq \g$, meaning that $N_{i,j}\in C_i\cap C_j$.
For example, if $N$ is an internal node of $Y$, then $\{ N \}$ is a cyclic set.
Observe that, if $Y$ is not of compact type, then $Y$ always admits
some cyclic set of nodes. Conversely, a curve of compact type does not admit
any cyclic set of nodes.

\proclaim Lemma \ecc.
Let $\Sigma$ be the complement in $Y_{sing}$ of a cyclic set of nodes.
Then $\Sigma$ is admissible and
$\es =g-\sum g_i-1.$
More generally, let $\Sigma$ be the complement in $Y_{sing}$ of a disjoint union of 
$n$ cyclic sets of
nodes. Then $\Sigma$ is admissible and
$\es \leq g-\sum g_i-n.$

\pf
If $\Sigma _1.....,\Sigma _n$ are cyclic sets of nodes of $Y$ such that $\Sigma _i\cap \Sigma
_j=\emptyset$ for all $i\neq j$, then their union contains an even number of external nodes on
every irreducible component of $Y$. Thus $\Sigma = Y_{sing} - \cup _1^n\Sigma _i$ is admissible
by definition. If $n=1$ so that $\Sigma$ is the complement of $\Sigma _1$ we
easily compute
$$
e_Y(\Sigma )= (\d - \# \Sigma _1)- (\g - \# \Sigma _1+1) + \g_c(Y) =\d -\g
+\g_c(Y) -1 
=g - \sum g_i -1.
$$
If $n $ is arbitrary, the inequality in the statement is straightforward. \qed

\

\df {\it Exponent set.}
 We shall call  the finite set of integers, that occur as exponents
of admissible sets of nodes,  the {\it exponent set} of $Y$, and we
will denote it by $E_Y$:
$$
E_Y:=\{e,\  {\rm{such\  that \  there\  exists }}\  \Sigma\in A_Y 
\ 
{\rm{\  with\  }}\es= e\}
$$

\ni
\rk
It is clear that $E_Y$ depends only on the combinatorial data of $Y$.
More precisely, if two curves $Y$ and $Y'$ have the same genus and the same dual graph,
then $E_Y=E_{Y'}$.

\proclaim Lemma \max .  
\item a) If $Y$ is  not of compact type, then $\{ g-\sum g_i -1, g-\sum g_i \}\subset
E_Y$.
\item b) If $\Sigma $ is admissible and 
$e _Y (\Sigma) = g -\sum g_i $, then $\Sigma = Y_{sing}$.
\item c)  If $N\in Y$ is an internal node, then either $g-\sum g_i -2 \in E_Y$
or $\Yn _N$ is of compact type and $E_Y =\{ 0,1\}$.

\pf
Since $Y_{sing}$ is admissible and $
e_Y(Y_{sing})=g-\sum g_i
$, we have that $g-\sum g_i\in E_Y$. 
To see that $g-\sum g_i-1\in E_Y$, observe that,
since $Y$ is not of compact type, $Y$  contains a cyclic set of nodes,
whose complement is admissible and has exponent precisely  $g-\sum g_i-1$,
by \ecc . This 
proves a).

For b) we have to show that if $\Sigma $ is admissible and its exponent is
 maximum  (that is $\es = e (Y_{sing})=g - \sum g_i$) then $\Sigma$
contains every node of
$Y$.
By contradiction, let $\Sigma \neq Y_{sing}$ and let $Z=\Yn _{\Sigma}$; the curve $Z$ is
singular. By Property D,
$\emptyset \in A_Z$ hence $Z$ is free from separating nodes,
thus $Z$ is not of compact type and  it contains a cyclic set of  nodes
$\Sigma _0$. We have thus a cyclic set of nodes 
$\Sigma _0\subset Y_{sing}$ such that
$\Sigma \subset \Sigma _0^c$.
Hence we have, by Property A,
$$
\es \leq e_Y(\Sigma _0^c)=g-\sum g_i -1
$$
which is a contradiction. 
To prove c) let $Z=\Yn _N$; since $N$ is internal, we have that 
$$
A_Y=\{ \Sigma , \Sigma \cup N, \  \forall \Sigma \in A_Z\};
$$
By Property B, $e_Y(\Sigma \cup N) = \es +1 = e_Z(\Sigma) +1$ hence 
$$
E_Y=\{ n , n+1, \  \forall n\in E_Z\}.
$$
If $Z$ is of compact type, $E_Z =\{ 0\}$ hence $E_Y =\{0,1\}$.
Otherwise, part (a) applied to $Z$ says that, $E_Z$ contains the number $n=(g-1)
-
\sum g_i -1$ (the arithmetic genus of $Z$ is, of course, $g-1$).
Hence $E_Y$ also contains $n= g-\sum g_i -2$ and we are done.
\qed

\

\ni{\bf \rel .}  Let $Y$ be stable. There is a relation between $E_Y$ and 
$S_Y$:
$$
L(S_Y)=\cases{\{ 2^n, \  \  \forall n \in E_Y,\  n\neq g - \sum g_i\} &  if \
$g_i=0\  \forall i$\cr
\{ 2^n, \  \  \forall n \in E_Y\} &  otherwise\cr}
$$

In fact,  $L(S_Y)$ is described in \lsy, and recall that for every $\xi \in S_Y$,
the support $\Sigma _\xi$ of $\xi$ is admissible. 
Now, (see [Co] section 6) if $\Sigma$ is admissible and $\Sigma \neq
Y_{sing}$, then there exist odd (and even) spin curves supported on $\Sigma$
(an equal number of odd and even, in fact).
By \max , $\es = g-\sum g_i$ iff $\Sigma = Y_{sing}$. A spin curve $\xi$ supported
on $Y_{sing}$ is given by a line bundle $L$ on $\Yn$ which restricts to a
theta-characteristic $L_i$ on every irreducible component $C_i^{\nu}$ of $\Yn$. For
$\xi$  to be odd it is necessary and sufficient that $L_i$ be an odd 
theta-characteristic on $C_i^{\nu}$ for an odd
number of components of $\Yn$. 
Therefore there exist odd spin curves supported on $Y_{sing}$, unless
 all components of $\Yn$ are rational;
in fact there are obviously no odd theta-characteristics on 
 $\P ^1$.
In other words, if $g_i=0$ for every $i$, a spin curve of exponent $g$ is
necessarily even.

\

\proclaim Corollary \ctp . A  curve $Y$ is of compact type if and
only if $E_Y =\{0\}$.

\rk By \rel, if $Y$ is stable, this is equivalent to saying that $Y$ is of
compact type if and only if $S_Y$ is reduced.

\pf If $Y$ is of compact type, the only admissible set is $Y_{sing}$
(by Property C),
whose exponent is $0$. The converse follows  from
 Lemma \max. \qed

\proclaim Lemma \inv . Let $Y$ be a connected, nodal  curve and let $X$ be its
stable model.
Then $E_X=E_Y$.

\pf
Let $\s : Y\la X$ be the natural map, contracting all smooth rational components
of $Y$ meeting their complementary curve in less than 3 points.
There is a natural bijection $\b$ between $A_Y$ and $A_X$
$$
\matrix{\b:& A_Y & \mapright {} & A_X\cr
\;&\Sigma &\mapsto &\s(\Sigma)\cap X_{sing} }
$$
whose inverse is, denoting by $Y_{sep}$ the subset of all separating nodes of
$Y$,
$$
\matrix{\b ^{-1}:& A_X & \mapright {} & A_Y\  \  \cr
\;&\Sigma ' &\mapsto &[\s ^{-1}(\Sigma ')\cap Y_{sing}]\cup  Y_{sep}}
$$
It is a trivial verification to show that the two above maps are each other
inverse and that they preserve the exponents, that is
$e_Y(\Sigma) = e_X(\b (\Sigma ))$.
\qed

\

\ni {\bf\esp .}
We now describe   $E_X$ for a split curve $X$.
 Of course $e(X_{sing}) = g$. Let  $\Sigma$ be a subset of nodes such that
$\Sigma
\neq X_{sing}$; then $e_X(\Sigma) =\#\Sigma$ and 
$\Sigma$ is admissible if and only if $\  \#\Sigma \not\equiv g \  $ mod  $2$.
Thus if $g$ is odd, the exponents appearing in $E_X$ are $0,2,4,....,g-1,g$
(in particular, the only odd exponent is $g$). Symmetrically for 
 $g$  even.
Summarizing:

{\it Let $X$ be a split curve of genus $g$. Then}
$$
E_X=\cases{\{ 0,2,4,....,g-3,g-1,g\} & {\it if \ $g$ is odd;}\cr
\{ 1,3,5,....,g-3,g-1,g\} & {\it if \  $g$ is  even.}\cr}
$$

\

\ni
We shall prove that the converse is also true, in other words 
 split curves are identified in $\Mg$
by their exponent set. 
Something more precise is true.

\proclaim Theorem \comb .  Let $Y$ be a DM stable curve of genus $g$. 
Assume that  $g\in E_Y$  and that $g-2\not\in E_Y$.
Then either $Y$ is a split curve, or  $Y$ is the ``polygonal" curve of genus $3$.

Where recall that the polygonal curve of genus $3$ is  the nodal curve made of four copies of
$\P ^1$ meeting pairwise in one point (having a total of six nodes). Its canonical
model is the union of four general lines in
$\P ^2$, whence the name ``polygonal".

 We shall apply the Theorem above to obtain the following
crucial result:

\proclaim Corollary \split .
Let $X$ and $Y$ be stable curves of genus $g$; assume that $X$ is split.
(a) If $E_X=E_Y$, then $Y$ is split. (b) If $L(S_X)=L(S_Y)$ then $Y$ is split.

\ni
{\it Proof of the Corollary.} We can exclude that $Y$  is the polygonal curve of genus $3$, 
in fact, in  that case
$E_Y=\{2,3\}$
whereas the exponent set of a split curve $X$ of genus $3$ is $\{0,2,3\}$.
Similarly, $L(S_Y)$ does not contain $1$, whereas $L(S_X)$ does.

Assume now that $E_X=E_Y =E$. If $X$ is split, then
$E$  satisfies 
the assumptions of  Theorem \comb (see above), hence $Y$ is
split.

This proves (a). 
Let us show that (b) implies (a). For any stable curve $Z$, the relation
between
$E_Z$ and $L(S_Z)$ has been explained at the beginning of \rel.
This yields that $E_X$ and $E_Y$ contain precisely the same integers $n$ such that $n\leq g-1$.
It remains to show that $g\in E_Y$. Notice that $g-2\not\in E_Y$ and $g-1\in E_Y$; therefore,
if $g\not\in E_Y$,  we get a
contradiction to
\max (a) ($Y$ is not of compact type by \ctp). \qed

\

\ni
{\it Proof of the Theorem.}
Recall the basic notation:
$Y=\cup _{i=1}^\g C_i$. The proof  is in three steps. 

\

\ni
{\it Step  1:  $Y$ is a union of smooth rational
components.} 
The fact that $g\in E_Y$  implies  that $e_Y(Y_{sing})=g$ (see  \max )
and that all irreducible components of $Y$ have geometric genus $0$.
If $Y$ contained an internal node, by  \max , part c),   $E_Y$  would
contain $g-2$, which is not the case.

\

\ni
{\it Step  2:  It suffices to assume that $Y$ is free
 from separating nodes.} 
Given our stable curve $Y$, denote by $Y^*$ a stable curve of genus $g$ obtained by smoothing
out all the separating nodes of $Y$. This is to say that if $N$ is a separating node, and
$N=C_1\cap C_2$, the curve $Y^*$ is obtained by replacing $C_1\cup C_2$ (which is a curve of
arithmetic genus $0$ by the Step 1) with a smooth rational curve $D$.

Notice that, unless $Y$ is free from separating nodes, $Y^*$ contains a separating component:
$D$ (that is
$D^c$ is disconnected). 

By Property C there is a natural bijection between the
admissible sets of $Y$ and those of
$Y^*$: to an admissible set $\Sigma$ of $Y$ (which must contain all the separating nodes of
$Y$) we associate the admissible set $\Sigma ^*=\Sigma -\{N:N \  \rm{is\ separating}\}$ on
$Y^*$. By Property B we have
$
\es = e(\Sigma ^*),
$
hence  $E_Y = E_{Y^*}$.

If we show that $Y^*$ is a split curve, then $Y^*$ has no separating
components, hence $Y=Y^*$ and we are done.
From now on we shall assume that $Y$ is free from separating nodes.
 
\

\ni
{\it Step  3.} 
We can assume that $Y$ is free from internal and separating nodes.
We proceed by induction on $g$.
The case $g=3$ is the first case to  be treated. 
A straightforward case by case analysis (there are only 3 cases) yields the result.
By the previous steps, there are only three cases.
If $Y$ is not split, then $Y$ has more than two components.
If $Y$ has three components, then its dual graph  is uniquely determined and so is $E_Y$
(compare with the Remark in \cyc);
$Y$ has $5$ nodes and
$E_Y=\{1,2,3\}$. If $Y$ has four components, then $Y$ is the polygonal
curve and $E_Y=\{2,3\}$.
There are no other cases.

Assume $g\geq 4$. Let $N$ be any node of $Y$ and let $Z=\Yn _N$ be the normalization of $Y$ at
$N$.
Since $N$ is  not an internal node of $Y$, to complete the proof of the Theorem
it suffices to prove that $Z$ is a split curve.

Notice that $Z$ is connected ($Y$ is free from separating nodes), nodal,
and its arithmetic genus is $g-1$.
Moreover, $Z$ is a union of smooth, rational components (because $Y$ is) and hence
$g-1 \in E_Z$.

Let now $\Sigma \in A_Z$, then clearly $\Sigma \cup N\in A_Y$ and 
$e_Z(\Sigma)=e_Y(\Sigma \cup N) -1$.
We conclude that
$E_Z\subset E_Y -1$, that is, if $n\in E_Z$ then $n+1\in E_Y$.
In particular, $g-3 \not\in E_Z$, because, by hypothesis, $g-2 \not\in E_Y$.

Let now $X$ be the stable model of $Z$
(possibly equal to $Z$). By Lemma \inv,
$E_X=E_Z$. Therefore $X$ satisfies the assumptions of the Theorem. By induction we obtain
that $X$ is a split curve of genus $g-1$, unless $X$ is the polygonal curve of genus $3$.

The rest of the proof consists in showing that $Z$ is stable (i.e. $X=Z$)
and that $Z$ is not the polygonal curve of genus $3$.

Assume first that $X$ is split.

By contradiction, assume that $Z$ is not stable; then $Z$ has one or two destabilizing
components (since $Z$ is the partial normalization of the stable curve $Y$ at  a unique
node
$N$). Hence $Z$, and likewise $Y$, has three or four irreducible components.
We need to distinguish the components of $Y$ from those of $Z$; let us denote by $C_i^Z$
the irreducible component of $Z$ which naturally 
(via the normalization at $N$) corresponds to $C_i$ in  $Y$.

We can write $X=D_1\cup D_2$ with $\#D_1\cap D_2 =g$. 
Let $C_1$ and $C_2$ be the two components of $Y$ that correspond to $D_1$ and
$D_2$. Then $\# C_1\cap C_2 \geq 2$; pick two nodes in $C_1\cap C_2$ and denote them by
$N_{1,2},N_{2,1}$. Obviously the set $\{N_{1,2},N_{2,1}\}$ is a cyclic set of nodes of $Y$.

The fact that $X$ is split implies that every
destabilizing component of $Z$ meets both $C_1^Z$ and $C_2^Z$ in one point,
thus (on $Y$) we have that, for all $i\geq 3$ (so that $C_i^Z$ is destabilizing in $Z$),
$C_i\cap C_1\neq \emptyset$ and $C_i\cap C_2\neq \emptyset$.
There are two possible cases, according to the number of destabilizing
components of $Z$, which is 1 or 2.

First case: $Z$ has one destabilizing component, $C_3^Z$.
In this case, (up to switching $C_1$ and $C_2$)  $N\in C_2\cap C_3$ so let us
rename $N=N_{2,3}$. Since $C_2^Z\cap C_3^Z\neq \emptyset$ there must be another node in
$C_2\cap C_3$,
let it denote by $N_{3,2}$. The set $\{N_{2,3},N_{3,2}\}$ is a cyclic set of
nodes of $Y$.

Let now $\Sigma$ be the complement of the two cyclic sets we constructed above, that is
$$
\Sigma :=Y_{sing}-\bigl( \{N_{1,2},N_{2,1}\}\cup\{N_{2,3},N_{3,2}\}\bigr).
$$
By Lemma \ecc, $\Sigma\in A_Y$ and its exponent is easily seen to be equal to $g-2$.
Then $g-2\in E_Y$, a contradiction.

In the second case $Z$ has two destabilizing components, $C_3^Z$ and
$C_4^Z$. Now $N\in C_3\cap C_4$, so we rename it: $N=N_{3,4}$.
Since $C_3^Z$ and $C_4^Z$ both meet $C_1^Z$ and $C_2^Z$, we find
(in $Y$) nodes $N_{2,3}\in C_2\cap C_3$
and
$N_{4,2}\in C_2\cap C_4$.
The set of nodes $\{N_{2,3},N_{3,4},N_{4,2}\}$ is cyclic. Continuing as in the previous case,
we let $\Sigma$ be the set
$$
\Sigma :=Y_{sing}-\bigl( \{N_{1,2},N_{2,1}\}\cup\{N_{2,3},N_{3,4},N_{4,2}\}\bigr).
$$
which is admissible, by \ecc, and whose exponent is $g-2$. A contradiction.

We conclude that $Z$ is stable, which is what we wanted.

To complete the proof we must show that $X$ is not the polygonal curve of genus
$3$.
This is done  just as above,  showing,
by a trivial case by case analysis, that $Y$ contains two disjoint cyclic sets of
nodes  whose complement has exponent $g-2$.
\qed

\

\cl{\bf{4. Split Curves}}

\

\ni
{\bf \ele.}
The theta-hypersurface of a
projective split curve  is described in details in [C], where it is proved (Theorem
5) that split curves are uniquely determined by their theta-hypersurface, among all
curves in $V$. In this section, we shall give a sharper version of such a result
(\AMS ).
Furthermore, in section \GA, we will  study the behaviour of certain configurations
of  theta-hyperplanes of a split curve, under the natural action of $G=PGL(g)$.
The results of this  analysis will  be applied in the sequel.

It is worth pointing out a  nice  feature of split curves which is to be used
often: 
the projection of a split curve in $\P^{g-1}$ from every subset of $i\leq g-3$ of its
nodes is a split curve in
$\P ^{g-1-i}$.

Let $H$ be a theta-hyperplane of  type
$i$ of a projective split curve $X$, thus $H$  contains exactly $i$ nodes of $X$
($0\leq i \leq g-1$) and it is tangent to $X$ at $g-1-i$ smooth points,
equally distributed among the two irreducible components of $X$. The set of all
hyperplanes of type $i$ is denoted by $\T _i(X)$ and its cardinality  by 
$
t_i (X) = \#\T _i(X)
$. 
It is easy to see that $\T _i(X)$ is empty for all $i$
having the same parity of $g$.

We need the following 

\proclaim Claim. For every projective split
curve $X$, we have
$
t_{g-3}(X)=\tg3 \  \  and\  \   t_{g-1}(X) = {g+1 \choose g-1}.
$

For the first formula,
it suffices to check that for any subset $\Sigma $ of $g-3$ nodes of $X$, there
are $4$ theta-hyperplanes of type $g-3$ containing $\Sigma$,
and containing no other node of $X$. If $g=3$,
we are simply saying that two plane conics meeting transversely have exactly 4
distinct tangent lines in common (which is clear by looking at the duals)
and that these 4 tangent lines do not go through the common points
(also easy, see [CS]). The
general case is done by projecting $X$ 
from
$\Sigma$ to a split curve  of genus $3$ in $\Pl$.

The second formula  is obvious: the nodes of $X$ are in general linear position,
therefore there is a unique  theta-hyperplane of type $g-1$ for every set of $g-1$
nodes of $X$.

\

\ni
{\bf\HYP .} The dualizing sheaf of a split curve $Y$ might fail to be
very ample: if this is the case,  the canonical image $W$ of $Y$ is a double
rational normal curve in $\P ^{g-1}$. 
In other words, $Y$ has a degree $2$ morphism onto
$\P ^1$ and it
 is thus in the closure, in
$\Mg$, of the locus of  smooth,
hyperelliptic curves (see [HM] 3.159 and 3.160). 
A split curve of this type will be called a {\it hyperelliptic split curve}.
Of course in such a case, the canonical model
$W$ is  not theta-generic and its theta-hypersurface is not
consistently defined.

We shall  leave the projective point of view, for a moment, and examine the situation
abstractly. Let $Y$ be an abstract split curve and let $S_Y$ the scheme
parametrizing its odd spin curves. 
We are interested in the spin curves parametrized by
$\hats _Y$, whose exponent is at least $g-3$ and whose
space of global sections has dimension $1$ (see \rig ).

Let $\xi\in \hats_Y$ have exponent  $g-3$; then $\xi$ is supported on a set of
$g-3$ nodes of
$Y$. Moreover, every set $\Sigma \subset Y_{sing}$ such that
$\#\Sigma =g-3$ is the support of exactly $4$ odd, distinct
spin curves (this follows from [Co]). Therefore the number of spin
curves of $Y$ having exponent $g-3$ is $\tg3$.

There are no spin curves having exponent $g-2$. The set of $\xi
\in S_Y$ having exponent $g-1$ is easily seen to have cardinality equal to the
number of subsets of
$g-1$  nodes of $Y$, that is ${g+1 \choose g-1}$.

Summarizing, the zero-dimensional scheme $\hats _Y$ has 
 $
\tg3  + {g+1 \choose g-1}$ irreducible components,
its length is equal to $2^{g-3}\cdot
\tg3  + 2^{g-1}{g+1 \choose g-1}$ and $L(\hats _Y ) =\{2^{g-3} , 2^{g-1}\}$.

Let $\sigma :Y\la W\subset \P^{g-1}$ be a canonical model of the split curve $Y$,
thus
$W$ is either a split curve
or $\sigma$ is a two-to-one morphism onto a   rational normal curve $C$. Let
$$\{N_1,....,N_{g+1}\}:=\s (Y_{sing})$$ 
so that
the points $\{N_1,....,N_{g+1}\}$ are
all distinct and, of course, in general linear position,
since they lie on a rational normal curve.

Let $\Y \la T$ be a general one-parameter deformation of $Y$ to smooth,
theta-generic curves, and let $\Y \la \W\subset \P ^{g-1}\times T $ be a
canonical image inducing
$\sigma$ on the central fibers.
Then   $\t ^{\W}(W)$ is defined as the limit
of the theta-hypersurfaces of the generic fibers (see \AP).

We shall now describe a configuration of hyperplanes $\tdeg\subset \t ^{\W}(W)$
which corresponds to $\hats _Y$ and is independent of the choice of $\Y$ and
$\W$.

We  start  by writing   explicitely
$$
\tdeg :=\{2^{g-3}\that ,2^{g-1}\thatt\}
$$
where $\that  $ will correspond to spin curves having exponent $g-3$ and 
$\thatt$ to spin curves having exponent $g-1$.

We have two cases: either $\s$ is an isomorphism, or it is not.
In the first case, $W$ is a projective split curve,  we shall 
simplify the notation, writing $\tdeg = \hatt (W)$, and we have 
$$
\tdeg =\hatt (W):=\{2^{g-3}\T_{g-3}(W),2^{g-1}\T_{g-1}(W)\}
$$
defined at the beginning of the section. The hyperplanes in
$\T_{g-3}(W)$ correspond naturally to spin curves of exponent $g-3$ and
similarly, theta-hyperplanes in 
$\T_{g-1}(W)$
correspond to spin curves of exponent $g-1$.

If  $K_Y$ is not very ample and the image of $\s$ is a
rational normal curve $C$,  we define
$$
\that :=  \{ <N_{i_1},....,N_{i_{g-3}},T_{N_{i_{g-2}}}C>,
\forall i_1,....,i_{g-2}\  \  {\rm  s.t.} \   \  i_j\neq i_ h \}
$$
where by $T_NC$ we denote the tangent line to $C$ at $N$. The other piece of $\tdeg
$ is defined by
$$
\thatt =\{ <N_{i_1}+...+N_{i_{g-1}}>\  {\rm s.t.}\  \  1\leq
i_1<....<i_{g-1}\leq g+1 \}
$$

\proclaim Lemma \HAT . Let $Y$ be a split curve and let $Y\la W$ be a canonical
model. The configuration
$\tdeg$ defined above has the following naturality property.
For every one-parameter deformation
$\Y \la T$ of $Y$ to  theta-generic curves, and for any choice of
a canonical model
$\Y \la \W \subset \P^{g-1}\times T$ restricting to the above  $Y\la W$,
we have that
$$
\tdeg \subset \t ^{\W}(W).
$$
Moreover, there is a natural bijection between the irreducible components of
$\tdeg$ and the irreducible components of $\hats _Y$.

\pf
Consider first a general one-parameter deformation $\Y$ of $Y$.
Then, by \cor (a),  we have a natural, birational $T$-morphism $\mu:S_{\Y}\la
J_{\W}$; denote by $\hat{\mu}:\hats _Y\la J_{\W}$ its restriction to $\hats _Y$.
We are stating  that its image $\hat{\mu}(\hats _Y)$ is independent of the choice
of the family, and it corresponds to $\tdeg$. 

If $W$ is a split curve, this follows from the discussion at the beginning of
the section. 

Otherwise  $Y$ is hyperelliptic, thus, as we mentioned above,
 $Y$ can be expressed 
as  the specialization of a family
of smooth hyperelliptic curves $Y_{\eta}$. We shall use this to describe $\hats _Y$.

The Weierstrass points $W_1,....W_{2g+2}$ of the generic, smooth
fiber
$Y_{\eta}$,  specialize in pairs to the $g+1$ nodes of $Y$.
Let us fix the notation so that, for $i=1,....,g+1$, the pair $(W_i,W_{i+g+1})$
specializes to the node that maps to $N_i$. Recall now that for smooth,
hyperelliptic curves, a complete description of the set of
theta-characteristics is given in [Mu2], Proposition 6.1.  
The hyperplanes   of
$\that$ correspond to the specialization of odd theta-characteristics 
of type ${\cal
O}_{Y_{\eta}}(W_{i_1}+W_{i_2}+...+W_{i_{g-3}}+W_{i_{g-2}}+W_{i_{g-2} +g+1})$,
with $1\leq i_1<....<i_{g-2}\leq g+1$.

By \rig, there is no ambiguity for the choice of the hyperplane corresponding
to a given $\xi \in \hats _Y$.

In  a completely similar fashion we  deal with  spin curves of exponent $g-1$.
This is in some sense easier, in fact such spin curves 
are supported on sets of $g-1$ nodes (all of such sets will occur).
Moreover, for a given set $\Sigma$ of $g-1$ nodes of $Y$ there exists a unique,
odd spin curve supported on $\Sigma$ (by [Co]). As before, viewing
$Y$ as specialization of smooth hyperelliptic curves, the hyperplanes of $\thatt$
correspond  to the
specializations of theta-characteristics of type
$
{\cal
O}_{Y_{\eta}}(W_{i_1}+W_{i_2}+...+W_{i_{g-2}}+W_{i_{g-1}})
$
with $1\leq i_1<....<i_{g-1}\leq g+1$.

So our statement holds for a general $\Y$. 
Now an arbitrary one-parameter deformation of $Y$ (satisfying the assumptions)
can be obtained as a specialization of general ones, for which the result
holds. A standard specialization argument (as in
\num) gives the result in general. \qed

\

\ni
{\bf \GA.} We need to study the behaviour of the above defined configurations of
hyperplanes, with respect to the  natural action of
$G=PGL(g)$ on 
$Sym ^m (\P ^{g-1})^*$.
The spaces of configurations of hyperplanes are fully understood from
the Geometric Invariant Theory point of view.
In [GIT] a criterion for stability is proved,
 which we need to recall. 
For any $\Omega \in Sym ^m (\P
^{g-1})^*$,  define, for every
$h=0,...., g-2$, $\mu _h(\Omega)$ as the maximum multiplicity of an
$h$-dimensional linear subspace $\Lambda \subset \P ^{g-1}$ as 
a subscheme of $\Omega$ (viewed now as a degree $m$-hypersurface in $\P^{g-1}$).
Thus, for example, $\mu _0(\Omega)$ is the maximum order of a singular point of
$\Omega$ and $\mu _{g-2}(\Omega)$ is the maximum multiplicity of a component 
(that is, a hyperplane) of
$\Omega$.
The criterion for stability is the following ([GIT] Proposition 4.3):

\proclaim GIT Stability Criterion. $\Omega $ is  stable if
and only if for every
$h=0,....,g-2$ we have
$
\mu _h(\Omega) < Max_h(m),
$
where
$
Max_h(m):=m{g-1-h \over g}.
$

We are going to apply it to the next result.

\proclaim Lemma \git .
Let $C\subset \P ^{g-1}$ be a rational normal curve and let $N_1,....,N_{g+1}$
be distinct points on
$C$.  
Let $\Omega \in Sym ^m (\P ^{g-1})^*$ be a configuration of hyperplanes 
of type  a),  b)  or c) below. Then $\Omega$ is GIT-stable 
with respect to the natural action
of $G$ on $Sym ^m (\P ^{g-1})^*$.

\item a)  $m = \tg3$ and 
$\Omega = \T_{g-3}(X)$ where $X\subset \P^{g-1}$ is a projective split curve.

\item b) $m = \tg3$ and  $\Omega = \{
<N_{i_1},....,N_{i_{g-3}},T_{N_{i_{g-2}}}C>,
\forall i_1,....,i_{g-2}\  \  {\rm  s.t.} \   \  i_j\neq i_ h
\}$

\item c) $m = {g+1 \choose g-1}$ and 
$\Omega = \{ <N_{i_1}, ...., N_{i_{g-1}}>, \
 \    1\leq i_1 < ....< i_{g-1}\leq g+1 \}$

\

\rk
If $\Omega$ is as in b), we have, of course, that $\Omega =\that$. Similarly, if
$\Omega$ is as in c), then $\Omega = \T_{g-1}(X)=\T_{g-1}(Y,W)$.

\pf
To prove the result for cases a) and b), notice
that
for every $h$ we have
$$
\mu _h(\T _{g-3}(X))\leq \mu _h(\that);
$$
therefore, if $\that$ satisfies the criterion, $\T _{g-3}(X)$ also does.
It is thus enough to deal with case b).
Let us compute $\mu _h(\that)$.
The $h$-dimensional linear subspaces contained in $\that$,  and having the
highest
 multiplicity are, clearly, those spanned by $h+1$ among the generating
points
$N_i$. Of course, the definition of $\that$ being symmetric with respect to
the $N_i$, varying the
$(h+1)$-t-uple does not change the multiplicity. Therefore we can pick a
specific set, say $\{N_1,....,N_{h+1}\}$, of $h+1$ nodes and have
$$
\mu _h(\that ) = \# \Bigl\lbrace H\subset \that : \{N_1,....,N_{h+1}\}
\subset H\Bigl\}.
$$
Now there are two types of hyperplanes $H$ contributing to $\mu _h(\that )$:
either $H$ does not contain the tangent line $T_{N_{i}}C$ for any
$i=1,....,h+1$, or it does.
The hyperplanes $H$ of the first type are
$$
{g+1-(h+1)\choose g-3 - (h+1)}4 = {g-h\choose 4}4
$$
where the binomial coefficient is the number of $(g-3)$-uples containing the
fixed set $\{N_1,....,N_{h+1}\}$, and the coefficient $4$ is there because for
every chosen $g-3$ as above, we have a choice of $4$  points where the
hyperplane is tangent to $C$.

The hyperplanes of the second type contain the tangent line to $C$ at one of
the $N_i$, for $i\leq h+1$. The total number of them is 
$$
{g+1-(h+1)\choose g-3 - h}(h+1) = {g-h\choose 3}(h+1)
$$
where the binomial 
coefficient is the number of $g-3$-uples containing a fixed subset of $h$
elements in  the   set $\{N_1,....,N_{h+1}\}$;
the factor $(h+1)$ represents  the choice of the  $N_i$ such
that
$H$ contains the tangent line to $C$ at $N_i$.

In conclusion
$$
\mu _h(\that )=
{g-h\choose 4}4+{g-h\choose
3}(h+1)=
{(g-h)(g-h-1)(g-h-2)(g-2)\over 3!}
$$
On the other hand, the strict upper bound allowed for stability by the criterion
of GIT  is
$$
Max_h(m)=m{g-1-h \over g}=4{g+1\choose g-3}{g-1-h \over g}=
{(g+1)(g-1)(g-h-1)(g-2)\over 3!}.
$$
Since $0\leq h$, it is clear that for every $h$
$$
\mu _h(\that )<Max_h(m)
$$
and therefore $\that$ and $\T _{g-3}(X)$ are GIT-stable.

Now we treat  case c), which is simpler.
We have for every $h$
$$
\mu _h(\Omega )=
{g+1-h-1\choose g-1-h+1}=
{g-h\choose 2} = {(g-h)(g-h-1)\over 2}.
$$
On the other hand, the (strict) upper bound given by the criterion above is,
for every $h$
$$
Max_h(m)=m{g-1-h\over g}={g+1 \choose 2}{g-1-h\over g}={(g+1)(g-1-h)\over 2}
$$
it is thus evident that, for every $h\geq 0$,
$\mu _h(\Omega )<Max_h(m)
$, and hence $\Omega$ is GIT-stable.
\qed

\

\rk
The interested reader can generalize the above computation to see that such
stability results are special cases of a more general phenomenon, about the
stability of configurations of hyperplanes spanned by sets of points 
and tangent lines of a
rational normal curve.

\

\proclaim Corollary \sta . Let $X$ be a projective split curve;
let $Y$ be an abstract split curve and $Y\la W$ a canonical model.
Consider the  configurations $\hatt(X)$ and $\tdeg$ 
with respect to the natural action of $G$.
(a) They are
GIT-stable. (b) If they are in the same $G$-orbit,
then $W$ is a projective split curve.

\pf 
The proof of (a) is a straightforward application of the
criterion already  used for the proof of the previous result. 
Let $m_1 = 4{g+1\choose g-3}$, $m_2 = {g+1\choose g-1}$
and
$
m = 2^{g-3}m_1+2^{g-1}m_2.
$
By definition, 
$$
Max _h(m) =2^{g-3} Max _h(m_1)+2^{g-1} Max _h(m_2)
$$
where $0\leq h\leq g-2$.
Let now $\Omega \in Sym ^m (\P ^{g-1})^*$ be one of the two configurations
in the statement.
Then we can write $\Omega =\{ 2^{g-3} \Omega _1,2^{g-1} \Omega _2 \}$
where $\Omega _1$ is 
either $\T _{g-3}(X)$ or $\that$, and hence 
a configuration of $m_1$ hyperplanes of
 type a) or b) in the previous Lemma \git \  a); similarly, $\Omega _2 $ is  a
configuration of 
$m_2$ hyperplanes of
 type c) in \git .
We have
$$
\mu _h(\Omega) = 2^{g-3} \mu _h(\Omega_1)+2^{g-1} \mu _h(\Omega_2) <2^{g-3} Max
_h(m_1)+2^{g-1} Max _h(m+2) = Max _h(m)
$$
(where the inequality comes from \sta ), hence we are done.

For (b) we must prove  that the $G$-orbits of $\hatt ( X)$ and $\tdeg$ are
different, if $W$ is not a split curve. This
follows from
the fact that the two, regarded as hypersurfaces, have  different
singularities. More precisely, by looking at the points $N_i$,  one sees that
$$
\mu _0(\that )>\mu _0(\T_{g-3} (X))
$$
(recall in fact that if $H\in \T_{g-3}(X)$, then $H$ is tangent to $X$ at two smooth
points). Since, as already noticed, $\mu _h(\thatt ) =\mu _h(\T_{g-1} (X))$ for every
$h\geq 0$,
we conclude that 
$$
\mu _0(\tdeg )>\mu _0(\hatt (X))
$$
and hence the two configurations cannot possibly be projectively equivalent.
\qed

\

\ni{\bf \TP .}
We conclude with a strengthening of Theorem 5 of [C] which will be used later.
The improvement consists  essentially in the fact that to recover the curves it
is enough to consider theta-hyperplanes of type $g-3$ and $g-1$,
rather than all of them.

\proclaim Proposition \AMS. Let $X$ and $X'$ be two   split curves in 
$\P^{g-1}$. 
If $\hatt (X) = \hatt (X')$, then $X=X'$.

\pf
The proof is  similar to the proof of Theorem 5 in [C].
The argument there is divided into two parts; the first part 
shows  that, if the two curves have the same theta-hyperplanes of
type $g-1$, then they have the same singularities. We can here use that part,
since, by definition, $\hatt (X) =\hatt (X')$ implies that $X$ and $X'$ have the
same theta-hyperplanes of type $g-1$ (which are, of course, those hyperplanes
having multiplicity $2^{g-1}$).

Denote $X_{sing}=X'_{sing}=\{ N_1,....,N_{g+1}\}$, $X=C_1\cup
C_2$ and  $X'=C'_1\cup C'_2$, with $C_i$ and $C_i'$ rational normal curves.

Now let, for $j=1,2$
$$
\Lambda _j =<N_j,N_3,N_4,....,N_{g-3},N_{g-2}>
$$
thus $\Lambda _1$ and $\Lambda _2$ are two linear subspaces of dimension $g-4$
intersecting in $<N_3,....,N_{g-2}>$.
Let, for $j=1,2$,
$\pi _j:\P^{g-1}-->\Pl$
be the projection from $\Lambda _j$ onto $\Pl$
and let $X_j =\pi _j(X)$ and $X'_j=\pi _j(X')$.
$X_j$ and $X_j'$ are two plane quartics with the same singularities, moreover
$X_j$ is split, because $X$ is split, notice in fact that $X_j $ 
(respectively, $X_j'$) is the normalization of
$X$ (respectively, of $X'$) at the nodes $N_j,N_3,....,N_{g-2}$.
Since $X_{sing} = X'_{sing}$ we have that $X_j$ and $X'_j$ have the same
singularities. There is a natural bijection between the theta-lines of type $0$
of
$X_j$ (respectively, of  $X_j'$)
and the theta-hyperplanes of type $g-3$ of $X$ (respectively, of $X'$) that
contain
$N_j,N_3,....,N_{g-2}$. By assumption, 
$X$ and $X'$ have the same theta-hyperplanes of type $g-3$, therefore
$X_j$ and $X'_j$ have the same
theta-lines of type $0$. We conclude (see [CS])
that $X_j =X'_j$ for $j=1,2$.

We have thus proved that projecting $X$ and $X'$  to $\Pl$ from the same
$g-4$-dimensional linear subspace spanned by $g-3$ of their (common) nodes,
one obtains the same split curve in $\Pl$.
It is clear that, up to re-naming the nodes, we can chose $\Lambda _1$ and
$\Lambda _2$ so that 
for $j=1,2$ we have that
$\pi _j(C_1) =\pi _j(C'_1) $ and, of course, $\pi _j(C_2) =\pi _j(C'_2) $.

Let $S_j$ be the cone over $C_1$ with vertex $\Lambda _j$,
that is
$$
S_j :=\bigcup _{P\in C_1}<\Lambda _j,P>,
$$
then the fact that $\pi _j (C_1)=\pi _j(C'_1)$ implies that 
$C'_1\subset S_j$ for both $j=1,2$.
 We shall now show that this implies that $C_1 = C_1'$.
It suffices to show that
that
$$
S_1\cap S_2 \subset C_1\cup < \Lambda _1,\Lambda _2>.
$$
By contradiction;
suppose that there is a point $Q$ such that $Q\in S_1 \cap S_2$ but $Q\not\in
C_1\cup <
\Lambda _1,\Lambda _2>$. Then the linear space $<Q,\Lambda _j>$  must intersect
$C_1$ in a point $P_j$; we thus find two new points $P_1$ and $P_2$ 
lying in the intersection of $C_1$ with the hyperplane $H:=<Q,\Lambda _1,\Lambda
_2>$. Then
$$
\{ N_1,N_2,....,N_{g-2},P_1,P_2 \} \subset C_1\cap H
$$
that is, $\deg C_1\cap H\geq g$, which is not possible. Therefore $C_1=C_1'$.

Repeating the argument for $C_2$ and $C_2'$ we conclude
that $X=X'$. 
\qed

\

\cl{\bf  5. The  local picture near split curves}

\

The purpose of this section is to analyze the differential of 
$\t: V \to Sym^{N_g}(\P^{g-1})^\ast$ at a point parametrizing a  
 split curve, and to prove (in \dif ) that it is injective.

\ni
{\bf  \ET .  Vector bundles on $\P ^1$: elementary transformations.}
We recall a few well known facts concerning vector
bundles on rational  curves.  Denote by $C=\P ^1$ and  by $\lambda$  the line
bundle of degree 1 on $C$.  Let $E$ be a
vector bundle on C.
An {\it elementary transformation}
of
$E$ is a vector bundle $E'$ such that there is an exact sequence
$$
0 \to E' \to E \to {\Bbb C}_x \to 0 \leqno (1).
$$
where ${\C}_x$ is a torsion sheaf supported on a point $x\in C$ with
fiber ${\C}$. We have
$$
{\rm rk}(E') = {\rm rk}(E); \ \ {\rm deg}(E') = {\rm deg}(E)-1
$$
Let us denote by $\P(E)$ the projective bundle over $C$
associated to $E$. The  exact sequence  (1) corresponds to a
point of 
$\P(E)$, contained in the fiber $\P(E)_x$, called {\it the center of the elementary
transformation}
$E'$.  

It is well known that the projective bundle
$\P(E')$ is obtained from $\P(E)$ by blowing up the center  and then blowing
down the proper transform $\tilde F$ of the fiber containing it (see [Ma]); 
$\P(E)$ is obtained from $\P(E')$ by the analogous  ``inverse" process with
center the image  of $\tilde F$.  
For example, if $E$ has rank $2$, we  write  $\P (E) \cong {\Bbb F}_n$ and $\P
(E')
\cong {\Bbb F}_{n'}$, with the usual notation ${\Bbb F}_n:=\P(\l^0 \oplus \l^n)$. 
Of course, in our situation,
$|n-n'|=1$. More precisely, we have that $n=n'-1$ if and only if the center $p\in
\P(E)$ belongs to the 
 $(-n)$-curve of $\P (E)$, if and only if the center $q\in \P (E')$ of
the inverse process does not belong to the $(-n-1)$-curve of $\P (E')$. 

 An exact
sequence of vector bundles
$$
0 \to F \to E \to L \to 0
$$
with ${\rm rk}L=1$, defines a section $\sigma$ of $\P(E)$; $\sigma$
contains   the center of $E'$ if and only if $F\subset E'$. In this case  we have an
exact sequence:
$$
0 \to F \to E' \to L\lambda^{-1} \to 0
$$
In particular, if $E=F\oplus L$ then $E' = F\oplus L\lambda^{-1}$.

Given an exact sequence 
$$
0 \to E_n \to E \to {\cal F} \to 0 
$$
where 
${\cal F}= {\C}_{x_1} \oplus \cdots\oplus{\C}_{x_{n}}$ is a torsion sheaf
supported on $n$ points of $C$,
we say that $E_n$ is obtained from $E$ by a sequence of $n$ elementary
transformations, centered at $n$ points of $\P (E)$. 

Conversely, a set $\Gamma$ of $n$ points of $\P(E)$, no two on the same fiber,
defines a vector bundle $E_\Gamma$ endowed with an exact sequence
$$
0 \to E_\Gamma \to E \to {\cal F} \to 0 
$$
where ${\cal F}$ is as above. Of course, $E_{\Gamma}$ is obtained applying the
elementary transformations corresponding to the points of $\Gamma$.

We need to consider the following special case. Let $r={\rm rk}E$ and let $\Gamma
'$ be a set of $r$ points, all contained in a fiber of $\P (E)$, and
in general position in such a fiber. Then one can consistently define $E_{\Gamma
'}=E\otimes
\lambda ^{-1}$, so that $\P(E) \cong \P (E_{\Gamma '})$.
Similarly, if $n=rq$ and $\Gamma _0$ is a union of disjoint $r$-uples of points, as
the $\Gamma '$ above, then $E_{\Gamma _0}=E\otimes \lambda ^{-q}$.

We shall need the following very simple Lemma.

\proclaim Lemma \USC. Let ${\cal E}\la T$ be a family of vector bundles over $\P
^1$. Let $E_t=\oplus \lambda ^{n_t^i}$ be the generic fiber and let 
$E_0=\oplus \lambda ^{n_0^i}$ be the special fiber.
 Then 
\item a) ${\rm max}_i \{n_t^i\}\leq {\rm max}_i \{n_0^i\}$ and
${\rm min}_i \{n_t^i\}\geq {\rm min}_i \{n_0^i\}$.
\item b)If $E_0$ is balanced (i.e. $n_0^i=n_0^j$ for all $i,j$) then $E_t$ is
 balanced and $E_t\cong E_0$.

\pf
Let $M_t:={\rm max}_i \{n_t^i\}$ and $M_0:={\rm max}_i \{n_0^i\}$. By contradiction;
if $M_t> M_0$, we can consider the family ${\cal E}'\la T$ of vector bundles,
obtained by tensoring the given family
$\cal E$ by $\lambda ^{-M_t}$. Then $E'_0$ has no global sections, being a sum of
line bundles of negative degree. On the other hand $E'_t$ has one summand equal to
$\lambda ^0$, thus it has  non-zero sections. This is a contradiction. The 
statement about the minimum is obtained in the same way,
working on the sequence of dual vector bundles.
Part b) follows immediately fom a).
\qed

\

\ni
{\bf  \NB . Normal bundles of split curves.} Let now $g \ge 3$ and let $C \subset \P^{g-1}$ be a rational
normal curve. 
Let us consider the normal bundle of $C$:
$$
\N_C := \N _{C/\P ^{g-1}}\cong \bigoplus^{g-2}  \lambda^{g+1}
$$
We have 
$
{\rm deg}\N_C=(g-2)(g+1),
$
$
h^0(\N_C) = (g-2)(g+2)
$
and
$
\P(\N_C) \cong C \times \P^{g-3}.
$
This isomorphism can be interpreted in terms of the  identifications:
$$
\P(\N_C) \cong \P(\N_C(-1)) = \{(x,H): x\in C,\  T_xC\subset H\} \subset C\times
(\P^{g-1})^*
$$
($T_xC$ is the tangent line to $C$ at $x$) coming from the
surjections:
$$
\oplus ^g\O_C \to T_{\P|C}(-1) \to \N_C(-1)
$$
More precisely, we have:
$
\P(\N_C) \cong C \times |\lambda^{g-3}|.
$
Thus a fiber $\P(\N_C)_x$, $x\in C$, can be identified with the
complete linear system $|\lambda^{g-3}|$ cut out by all hyperplanes  containing
$T_xC$.

Let $X = C_1 \cup C_2 \subset \P^{g-1}$  be a   split curve
of genus $g$ and  let
$C_1 \cap C_2 = \{N_1,\ldots, N_{g+1}\} = X_{sing}$.
Let $\Omega = \T_{g-3}(X)$ be the set of all theta-hyperplanes of $X$ passing
through exactly $g-3$ nodes of $X$
(described in the previous section).
We write $\Omega = \Omega_1 \cup \ldots \cup \Omega_{{g+1\choose g-3}}$, where
$\Omega_i$ is the set of four
theta-hyperplanes  containing  a fixed subset $\Sigma _i$ of $g-3$ nodes. We
denote by $H_{i,j}$, $j=1,\ldots,4$,  the four hyperplanes of
$\Omega_i$.

Consider the everywhere nonreduced  zero-dimensional scheme 
$$Z:= [\bigcup_{H\in\Omega}H]\cap X_{reg}\subset X$$
and   let, for $k=1,2$,
$
Z_k = Z \cap C_k.
$
Then   $Z_k \subset C_k$  is supported at the
 $4{g+1\choose g-3}$ points  below
$$\Biggl\{ p_{i,1}^k,\ldots,p_{i,4}^k,\  i=1,\ldots ,\tg3\Biggr\}
$$
For a general split curve, such points are  distinct and $p_{i,j}^k\neq N_h$
for all choices of indeces.
We have natural surjections $\N_X\to T^1_Z$ and
$\N_{C_k} \to T^1_{Z_k}$, where $T^1_-$ denotes the first cotangent sheaf of $-$,
and we define
$$
\N'_X := \ker \{\N_X \to T^1_Z\}
$$
$$
\N'_{C_k} = \ker\{\N_{C_k} \to T^1_{Z_k}\}
$$
We have an  exact and commutative diagram for 
$k=1,2$:
$$
\matrix{
&&0&&0&&&& \cr
&&\downarrow&&\downarrow&&&& \cr
0&\to&\N'_{C_k}&\to& \N_{C_k}&\to&T^1_{Z_k}&\to&0 \cr
&&\downarrow&&\downarrow&&\Vert&& \cr
0&\to&(\N'_{X})_{|C}&\to&
(\N_{X})_{|C_k}&\to& T^1_{Z_k}&\to& 0 \cr
&&\downarrow&&\downarrow&&&& \cr
&&{\cal T}&=&{\cal T}&&&& \cr
&&\downarrow&&\downarrow&&&& \cr
&&0&&0&&&& } \leqno (3)
$$
where ${\cal T}= {\C}_{N_1} \oplus \cdots\oplus{\C}_{N_{g+1}}$ 
is a torsion sheaf.

The proof of Theorem \dif \  is based on the following key result

\proclaim Propositon \van . With the above notation, $H^0(X,\N'_X)=0$ for all  $g\geq
4$. 

\pf
We have an exact sequence:
$$
0 \to \N'_X \to (\N'_{X})_{|C_1}\oplus (\N'_{X})_{|C_2} 
\to Q \to 0\leqno (2)
$$
where  $Q = {\C}^{g-2}_{N_1} \oplus \cdots\oplus{\C}^{g-2}_{N_{g+1}}$,
it therefore suffices to prove that, for $k=1,2 $,
$$
H^0((\N'_{X})_{|C_k}) = 0 .
$$
To prove this, we will consider the above diagram (3).
From its first row, we see that $\N'_{C_k}$ is obtained
from
$\N_{C_k}$
 applying a sequence of $4{g+1\choose
g-3}$ elementary transformations, centered at the points of the set
$$\Gamma :=\{(p_{i,j}^k,\Sigma _i),\   i=1,....,{g+1\choose g-3},\  j=1,....,4\}$$

Let us fix now $k=1$ and drop the index $k$ for simplicity, denoting $C_1=C$.
Let $E:={\cal N}_C =\oplus ^{g-2}\lambda ^{g+1}$.
Let $n = 4{g+1\choose
g-3}$.

We can specialize our curve $X=C\cup C_2$ to a hyperelliptic split curve
(see \HYP), whose
canonical model $X_0$ is supported on $C$. 
We can do that by keeping the component
$C$ fixed (i.e. $C_2$ specializes to $C$) and by mantaining
the nodes  of every fiber
at  the points
$N_1,....,N_{g+1}$. 
Call $X_t$ the generic fiber and $X_0$ the special fiber.
Let $\Gamma _t$ be the $n$-uple of points in $\P(E)$
defined as the centers of the $n$ elementary transformation of the first
row of diagram (3), which we can now re-write for $X_t$
$$
0\la E_{\Gamma _t}\la E \la  T^1_t \la 0
$$
As $X_t$ specializes to $X_0$, 
the configuration of theta-hyperplanes $\T _{g-3}(X)$ specializes to a
configuration that has been defined and studied in the previous section. Namely,
call $Y_0$ the abstract, hyperelliptic split curve whose canonical model is
$X_0$. Then the limit configuration is $\Theta_{g-3}(Y_0,X_0)$
defined in \HYP, where $W=X_0$.
Therefore
 the set $\Gamma _t$ specializes to a set $\Gamma
_0\subset \P(E)$ which  has implicitely been
described in \HYP:
$$
\Gamma _0 = \{(N_j, \Sigma _i): \  \  \#\Sigma _i = g-3   ,\  \  N_j\not\in \Sigma
_i\}
$$
In particular, $\Gamma _0$ is entirely contained in the $g+1$ fibers over $N_1,....,
N_{g+1}$; each such a fiber contains  exactly ${g \choose g-3}$ 
 points  of $\Gamma _0$, and such points are in general position in every fiber
(because $C$ is a rational normal curve).

At this point, we need to distinguish three cases; first: $g\not\equiv 2 \  
({\rm mod}\ 3)$ and
$g\neq 4$; second:
$g=3x+2$ and third: $g=4$,

$\bullet $
Suppose  that $g\not\equiv 2 \  ({\rm mod}\ 3)$.

Then we can factor
$n=(g-2)(g+1) {g(g-1)\over 2}$, so that $\Gamma _0$ is a union of 
disjoint $g-2$-uples contained in a fiber of $\P(E)\la C$. We are thus in the
situation
described in \ET; we obtain that
$$
E_{\Gamma _0} = E\otimes \lambda ^{-(g+1) {g(g-1)\over 2}}=
\bigoplus ^{g-2}\lambda ^{a}
$$
where $a= g+1 - (g+1) {g(g-1)\over 2}$.
Thus $E_t$ specializes to the 
``balanced" vector bundle $E_0 =\oplus \lambda ^a$.
By \USC, $E_t$ is also balanced and isomorphic to $E_0$.

We have therefore proved that
$
\N '_{C}= \oplus ^{g-2}\lambda ^a.
$

A straightforward computation shows that, if $g\geq 5$, then $a<-g-1$

The exact sequence in the first column of (3) expresses $\N '_{C}$ as
obtained from $(\N'_{X})_{|C}$ by a sequence of $g+1$ elementary
transformations.
In particular, 
$$\deg (\N'_{X})_{|C} =\deg \N '_{C}+g+1.$$
If $g\geq 5$,
 we immediately deduce that each line bundle summand in the splitting
of 
$(\N'_X)_{|C}$ has negative degree, therefore
$H^0((\N'_X)_{|C})=0$ and we are done.

$\bullet$ Let $g=3x+2$.
Here the problem is that for every $h=1,....,g+1$ the number $f$ of points of $\Gamma
_0$ contained in the fiber over $N_h$ is not a multiple of $g-2$.
More precisely we have
$$
f:={g\choose g-3} = x + (g-2){3x(x+1)\over 2}
$$
Pick on each of these $g+1$ fibers a subset of $f-x$ points of $\Gamma _0$,
and call $\Gamma _0'$ the union of these $g+1$ subsets. Thus
$\Gamma _0'\subset \Gamma _0$,
$$\# \Gamma _0' = (g+1)(g-2){3x(x+1)\over 2}$$
and, as before,
$$
E_{\Gamma _0'}=E\otimes \lambda ^{ - (g+1){3x(x+1)\over 2}}=\bigoplus ^{g-2}\lambda
^{a'}
$$
where $a' = g+1 -(g+1){3x(x+1)\over 2}$.
As in the previous case, we  conclude that
$
\N '_{C}\subset  \oplus ^{g-2}\lambda ^{a'}.
$

One easily sees that, if $x\geq 1$, then $a'<-g-1$.
So we conclude  by looking at the first column of diagram (3),
exactly as we
did before.

$\bullet$ 
If $g=4$ the argument of the first case yields 
$
\N'_{C}= \l ^{-5}\oplus \l ^{-5}
$, which is not enough; we use a slightly different strategy.
We can write $ (\N'_X)_{|C}=\lambda ^u\oplus\lambda ^v$ with $u\leq v$ and $u+v =
-5$. 
To show that  $(\N'_X)_{|C}$ has no sections it suffices to show that $v-u< 5$.

Consider  the second column of diagram (3)
$$
0 \to\l ^{5}\oplus \l ^{5} \to (\N_{X})_{|C} = \lambda ^c\oplus \lambda ^d \to {\cal
T}
\to 0
$$
with $c+d = 15$ and $c\leq d$,
For $i=1,\ldots ,5$ let  $H_i$ be the plane spanned by the tangent lines at $N_i$ to
$C_1=C$ and $C_2$:
$$
H_i:=<T_{N_i}C_1,T_{N_i}C_2>.
$$
It is a local computation
 to check that ${\cal T}_{N_i}$  is
defined by the plane $H_i$.

Therefore the 5 centers of the elementary transformation above
are the points $(q_i,N_i)$ where $q_i\in C\cap H_i -\{N_i\}$ is the residual point
of intersection of $H_i$ with $C_i$.
It is easy to see that $p_i\neq p_j$.
We deduce that such centers
are not all contained in the same section.
In particular, $\P ((\N_{X})_{|C})\neq {\Bbb F}_5$,
hence $d-c < 5$.

Now we look at the second row of diagram (3).
By the same argument as in the previous part of the proof, using \USC a),
we can conclude that $v-u\leq d-c <5$; hence we are done.
 \qed

We can now prove the

\proclaim Theorem \dif . Let $X\subset \P^{g-1}$ be a general, projective split
curve. Then 
$\t$ is an immersion at (the point parametrizing) $X$.

\ni
{\it Proof.}
Let $X=C_1\cup C_2\subset \P^{g-1}$ be a general  split curve. 
To prove that $\t$ is an immersion at $X$ amounts to showing that 
$T_X\t^{-1}(\t(X)) = (0)$. Note that $\t^{-1}(\t(X))$ is a closed 
subfamily of $V$ consisting of canonical curves which
have $\t(X)$ as theta-hypersurface. Consider the subsheaf 
$\N_X'\subset \N_X$ studied above.
 Its sections define a subspace of $H^0(X,\N_X)=T_XV$ 
 consisting of first order deformations of $X$ which remain tangent 
to the hyperplanes of $\Omega$ at the points of $Z$. 
Therefore   
$$T_X\t^{-1}(\t(X)) \subset H^0(X,\N'_X)$$
 But 
$   H^0(X,\N'_X) = 0$ by \van, thus our statement is proved.
\qed

\

\

\cl{\bf 6. Characterizing canonical curves by their theta-hyperplanes}

\

For the reader's convenience, we restate our main Theorem, before proving it.

\proclaim Theorem \main .
Let $X$ and $X'$ be general canonical curves in $\P ^{g-1}$
having genus $g\geq  4$. If $\t
(X) =
\t(X')$ then $X=X'$.

\pf
By contradiction. If the statement is false we can find two families of
canonical curves as follows. The first family 
$$
\matrix{\X & \hookrightarrow {} & T\times \P^{g-1}\cr
\   \lmapdown{ } & &\cr
T& &}
$$
is a general one-parameter deformation of a general split curve $X$, 
with  $X_t$   smooth for $t\neq t_0$.  For every $t\in T$, we assume
that
$\psi _{\X}(t)\in V$.

The second family
$$
\matrix{\X '& \hookrightarrow {} & T\times \P^{g-1}\cr
\   \lmapdown{ } & &\cr
T& &}
$$
is such that the generic fiber $X'_t$ is a smooth, theta-generic canonical
curve,
 having the
same theta-hypersurface as $X_t$: $\t(X_t)=\t(X'_t)$  for every $t\neq t_0$.

We cannot, for the moment, say much about the special fiber $X'$ 
of the second family,
we shall use the existence of stable reduction to analyze it.
Modulo replacing $T$ by a finite covering, we can assume that $\X' \la
T$ admits stable reduction over $T$; let it denote by $\Y \la T$. We have that,
for
$t\neq t_0$, $X'_t$ is a canonical (isomorphic) model of $Y_t$.
The central fiber $Y$ of $\Y \la T$ is a  stable curve, and we have a
birational map
$\Y --> \X '$ which is an isomorphism away from the central fibers.

$\bullet$ We claim
that {\it $Y$ is a split curve.}

To prove that, consider first $J_{\X}\la T$ defined at the end of \BS,
and the family of odd spin curves $S_{\X} \la T$
 defined in \CS \  (abusing notation denoting by the same
symbol $
\X$ the abstract family and its canonical model). We can compare the two:
 by our genericity assumption, $ S_{\X} $ is a
smooth curve (by \exp ) dominating $J_{\X}$ (by \cor (a)).
Notice also that, by construction,
$
J_{\X}=J_{\X '};
$
in conclusion, we have a birational morphism of $T$-schemes
$$
S_{\X}  \la J_{\X }=J_{\X '}
$$
and by \cor (b),
$$
L(S_{X})= L(J_{X}) \eqno (1)
$$
Consider now the pull-back to $T$ of the space of odd spin curves of the
family $\Y\la T$, as usual denoted by $S_{\Y} $.
Away from their special fibers, $S_{\Y} $ and $J_{\X'}$ are isomorphic over
$T$; by what we have observed above, we have a birational $T$-morphism (recall
that
$S_{\X}$ is a smooth curve)
$$
\a : S_{\X}  \la S_{\Y} .
$$
Looking at the central fibers we get that
$$
L(S_{X})\geq L(S_{Y}),
$$
 equivalently, by (1)
$$
L(J_{X})\geq L(S_{Y}) \eqno (2)
$$
Now we apply  \num , with $\W = \X '$.
We get that
$$
L(S_{Y})\geq  L(J^{\X'}_{X'})
$$
but $J_{\X'}=J_{\X}$ and hence, of course, $J_{X}=J^{\X'}_{X'}$ so that we
get
$$
L(S_{Y})\geq  L(J_{X}) \eqno (3)
$$
Combining (2) and (3) (see \PO) and using (1), we conclude that
$$
L(S_{X})=L(S_{Y}).
$$
Finally, $X$ being a split curve, we use \split (b)
\  to conclude that $Y$ is
a split curve. The claim is thus proved.

Let now $Y\la W$ be a canonical map. 
Let 
$$
\Y \la \W \hookrightarrow {}  T\times \P^{g-1}
$$
be a canonical image of $\Y \la T$, whose restriction
to the central fibers is the above
$Y\la W$.    If $t\neq t_0$, $W_t$ is an isomorphic
model of
$Y_t$ and it is thus projectively equivalent to $X_t'$. That is, there exists
a morphism $\gamma :T - \{t_0\} \la G$ 
such that 
$$X'_t = W_t^{\gamma (t)}.$$
Hence
$$ \t (X'_t) =\t (X_t)= \t (W_t)^{\gamma (t)}.$$

By \HAT, we can consider the distinguished subconfiguration of
theta-hyperplanes 
$\hatt (X)\subset
\t (X)$; by definition,
$\hatt(X)$ is the subconfiguration (with multiplicities) of all components of
$\t (X)$ having multiplicity at least
$2^{g-3}$. Therefore it is
 the
specialization of a well defined subconfiguration of hyperplanes $\hatt_t\subset
\t (X _t)$. Similarly, consider the distinguished configuration
$\tdeg \subset \t^{\W}(W)$. By construction,  we have that
$\tdeg$ is the specialization of the family of configurations
$
(\hatt_t)^{\gamma (t)}.
$
Recall that, by \sta , $\hatt (X)$ and $\tdeg$ are GIT-stable points in $Sym
^m(\P^{g-1})^*$. We have just seen that they are specializations of two
$G$-conjugate families.
Therefore (since they are GIT-stable)
they are themselves $G$-conjugate, that is,
there exists an element
$g_0\in G$ such that
$$
\hatt (X)=(\tdeg )^{g_0}.
$$
Now this implies that $W$   is a split curve, by
\sta , (b). Hence
 $\tdeg = \hatt(W)$ and
$$
\hatt (X)=\hatt(W)^{g_0}=\hatt(W^{g_0}).
$$
Summarizing, $X$ and $W^{g_0}$ are split curves such that
$
\hatt (X)=\hatt (W^{g_0}).
$
We are thus in the position of applying \AMS , to conclude that 
$$
W^{g_0} = X.
$$
Now $W^{g_0}$ is the central fiber of the family over $T$ whose generic fiber
is 
$W_t^{\gamma (t)}$ . By Theorem \dif, we get that $W_t^{\gamma (t)}= X_t$,
hence we are done, since by construction $W_t^{\gamma (t)}= X'_t$.
\qed

\

\

\ni
{\bf References }

\smallskip

\ni
[ACGH] E. Arbarello, M. Cornalba,  P. Griffiths, J. Harris: 
{\it Geometry of Algebraic Curves 1} 
Grundlehren der Mathematischen Wissenschaften, 267. Springer, New York-Berlin, 1985

\smallskip

\ni
[CS] L. Caporaso, E. Sernesi: {\it Recovering plane curves from their bitangents.}
Preprint Alg-Geom AG/0008239 To appear in Journ. of Alg. Geom.

\smallskip

\ni
[C] L. Caporaso: {\it On modular properties of odd theta characteristics}
Advances in mathematics motivated by physics, pp.101-114, Contemporary
Mathematics series AMS 

\smallskip

\ni
[Co] M. Cornalba: {\it Moduli of curves and theta-characteristics.}
Lectures on Riemann surfaces (Trieste, 1987), 560--589,
World Sci. Publishing, Teaneck, NJ, 1989

\smallskip

\ni
[HM] J. Harris, I. Morrison: {\it Moduli of curves.} Graduate texts in Math.
187, Springer (1998)

\smallskip

\ni
[L] D. Lehavi: {\it Any smooth plane quartic can be reconstructed from its
bitangents.} Preprint AG/0111017

\smallskip

\ni
[Ma] M. Maruyama : {\it Elementary transformations in the theory of algebraic
vector bundles.} Algebraic geometry (La Rabida, 1981), 241-266, Lecture Notes in
Math. 961, Springer 1982. 

\smallskip

\ni
[Mu1] D. Mumford: {\it Curves and thir Jacobians.}
The University of Michigan Press, Ann Arbor,  1975.

\smallskip

\ni
[Mu2] D. Mumford: {\it Tata lectures on theta II.}
Progress in mathematics 43, Birkh\"auser 1984.

\smallskip

\ni
[Mu3] D. Mumford: {\it Tata lectures on theta III.}
Progress in mathematics 97, Birkh\"auser 1991. 

\smallskip

\ni
[GIT] D. Mumford, J. Fogarty, F. Kirwan:   {\it Geometric Invariant Theory.}
(Third edition) E.M.G. 34 Springer 1994.

\smallskip

\ni
[S] R. Salvati Manni: {\it On the differential of applications defined  on moduli
spaces of p.p.a.v. with level theta structure.}  Math. Z. 221 (1996),
231-241.

\end